\RequirePackage{fix-cm}
\documentclass[smallextended]{svjour3}

\usepackage[utf8]{inputenc} 
\usepackage[T1]{fontenc}    
\usepackage{hyperref}       
\usepackage{url}            
\usepackage{booktabs}       
\usepackage{amsfonts}       
\usepackage{nicefrac}       
\usepackage{microtype}      
\usepackage{lipsum}

\usepackage{mathptmx}
\usepackage[misc]{ifsym}

\usepackage{graphicx,mwe}
\usepackage{amsmath}
\usepackage{hyperref}
\usepackage{subcaption}
\usepackage{placeins}
\usepackage{bm}
\usepackage[]{natbib}

\usepackage{pdflscape,rotating}

\usepackage{algorithmic, algorithm2e}
\RestyleAlgo{boxruled}
\usepackage{pifont}

\hypersetup{colorlinks,breaklinks,
           linkcolor=blue,urlcolor=blue,
           anchorcolor=blue,citecolor=blue}

\newcommand{\hb}{\mathbf{h}}

\newcommand{\ngg}{n_{\mathbf{q}}}
\newcommand{\np}{n_{\mathbf{\ppm}}}

\newcommand{\ppm}{{\bm\theta}}

\newcommand{\uu}{\mathbf{u}}

\newcommand{\Ym}{\pmb{\xi}}

\newcommand{\Exp}{\mathbb{E}}

\newcommand{\betaa}{\bm{\beta}}

\newcommand{\thetaa}{\bm\theta}

\newcommand{\xii}{\bm\xi}

\DeclareMathAlphabet{\mathcal}{OMS}{cmsy}{m}{n}

\usepackage{pgfplots}
\usepackage{mathptmx} 
\usepackage{multirow}

\newcommand\Tstrut{\rule{0pt}{2.6ex}}         
\newcommand\Bstrut{\rule[-0.9ex]{0pt}{0pt}}   

\usepackage{tikz}
\usetikzlibrary{matrix} 
\usetikzlibrary{arrows} 
\usetikzlibrary{arrows.meta}
\usetikzlibrary{calc} 
\usetikzlibrary{shapes}
\usetikzlibrary{fit}
\usetikzlibrary{positioning}
\usetikzlibrary{intersections,patterns,pgfplots.fillbetween}
\usetikzlibrary{calc}
\usetikzlibrary{decorations.pathreplacing,perspective}
\usepackage{tikz-3dplot}

\usepackage{ifpdf}
\ifpdf
  \usepackage{epstopdf}
\fi

\pgfmathdeclarefunction{gauss}{2}{%
  \pgfmathparse{1/(#2*sqrt(2*pi))*exp(-((x-#1)^2)/(2*#2^2))}%
}

\tikzstyle{block} = [draw,rectangle,thick,minimum height=2em,minimum width=2em]
\tikzstyle{sum} = [draw,circle,inner sep=0mm,minimum size=2mm]
\tikzstyle{connector} = [->,thick]
\tikzstyle{line} = [thick]
\tikzstyle{branch} = [circle,inner sep=0pt,minimum size=1mm,fill=black,draw=black]
\tikzstyle{guide} = []
\pgfplotsset{compat=1.16}

\pgfkeys{
  /pgf/arrow keys/.cd,
  pitch/.code={%
    \pgfmathsetmacro\pgfarrowpitch{#1}
    \pgfmathsetmacro\pgfarrowsinpitch{abs(sin(\pgfarrowpitch))}
    \pgfmathsetmacro\pgfarrowcospitch{abs(cos(\pgfarrowpitch))}
  },
}

\pgfdeclarearrow{
  name = Cone,
  defaults = {       
    length     = +3.6pt +5.4,
    width'     = +0pt +0.5,
    line width = +0pt 1 1,
    pitch      = +0, 
  },
  cache = false,     
  setup code = {},   
  drawing code = {   
    \pgfmathsetmacro\pgfarrowhalfwidth{.5\pgfarrowwidth}
    \pgfmathsetmacro\pgfarrowhalfwidthsin{\pgfarrowhalfwidth*\pgfarrowsinpitch}
    \pgfpathellipse{\pgfpointorigin}{\pgfqpoint{\pgfarrowhalfwidthsin pt}{0pt}}{\pgfqpoint{0pt}{\pgfarrowhalfwidth pt}}
    \pgfusepath{fill}
    \pgfmathsetmacro\pgfarrowlengthcos{\pgfarrowlength*\pgfarrowcospitch}
    \pgfmathparse{\pgfarrowlengthcos>\pgfarrowhalfwidthsin}
    \ifnum\pgfmathresult=1
      \pgfmathsetmacro\pgfarrowlengthtemp{\pgfarrowhalfwidthsin*\pgfarrowhalfwidthsin/\pgfarrowlengthcos}
      \pgfmathsetmacro\pgfarrowwidthtemp{\pgfarrowhalfwidth/\pgfarrowlengthcos*sqrt(\pgfarrowlengthcos*\pgfarrowlengthcos-\pgfarrowhalfwidthsin*\pgfarrowhalfwidthsin)}
      \pgfpathmoveto{\pgfqpoint{\pgfarrowlengthcos pt}{0pt}}
      \pgfpathlineto{\pgfqpoint{\pgfarrowlengthtemp pt}{ \pgfarrowwidthtemp pt}}
      \pgfpathlineto{\pgfqpoint{\pgfarrowlengthtemp pt}{-\pgfarrowwidthtemp pt}}
      \pgfpathclose
      \pgfusepath{fill}
    \fi
    \pgfpathmoveto{\pgfpointorigin}
  }
}

\title{Reliability-based Topology Optimization using Stochastic Gradients}

\tikzset{
  saveuse path/.code 2 args={
    \pgfkeysalso{#1/.style={insert path={#2}}}%
    \global\expandafter\let\csname pgfk@\pgfkeyscurrentpath/.@cmd\expandafter\endcsname
                           \csname pgfk@\pgfkeyscurrentpath/.@cmd\endcsname
    \pgfkeysalso{#1}},
  /pgf/math set seed/.code=\pgfmathsetseed{#1}}

\begin{document} 

\author{Subhayan De \and  Kurt Maute \and Alireza Doostan 
}


\institute{
	S. De \at
	Smead Aerospace Engineering Sciences Department, University of Colorado, Boulder, CO 80303, USA\\
	\email{Subhayan.De@colorado.edu}         
\and
 K. Maute\\ \email{Kurt.Maute@colorado.edu}\\ 
 A. Doostan (\Letter)\\ \email{Alireza.Doostan@colorado.edu}
}

\date{Received: date / Accepted: date}

\maketitle

\noindent\textbf{Dedication:} {This work is in memoriam of Raphael (Rafi) T. Haftka (1944 -- 2020), a pioneer of modern design optimization and design under uncertainty. His work has laid the foundation for today's and future research on design optimization, including the work presented in this paper. 
}

\begin{abstract} 
This paper addresses the computational challenges in reliability-based topology optimization (RBTO) of structures associated with  the estimation of statistics of the objective and constraints using standard sampling methods. The aim is to overcome the accuracy issues of traditional methods that rely on approximating the limit state function. Herein, we present a stochastic gradient-based approach, where we estimate the probability of failure at every few iterations using an efficient sampling strategy. To estimate the gradients of the failure probability with respect to the design parameters, we apply Bayes' rule wherein we assume a parametric exponential model for the probability density function of the design parameters conditioned on the failure. The design parameters and the parameters of this probability density function are updated using a stochastic gradient descent approach requiring only a small, e.g., $\mathcal{O}(1)$, number of random samples per iteration, thus leading to considerable reduction of the computational cost as compared to standard RBTO techniques. We illustrate the proposed approach with a benchmark example that has an analytical solution as well as two widely used problems in structural topology optimization. These examples illustrate the efficacy of the approach in producing reliable designs.

\end{abstract}

\keywords{Topology optimization \and Reliability estimation \and Stochastic gradients \and Optimization under uncertainty}

%


\section{Introduction}


When a designed structure is built and used in real-life applications, the loading, geometry, and material properties are typically different from the values used during the design process due to uncertainty. 
To achieve a robust design, i.e., a design whose performance does not depend {\it much} on stochastic variations, these uncertainties must be incorporated in the optimization process. For example, in design under uncertainty, the mean of the cost function is minimized subjected to constraints that are satisfied in expectation \citep{nikolaidis2004comparison,beyer2007robust,de2017efficient,Diwekar2020}. To reduce variability in the design's performance, a standard deviation or variance term can also be added to the objective \citep{beyer2007robust,dunning2013robust,de2020topology}. 
While the robust design procedure addresses the presence of uncertainty, the designer may be interested in limiting the failure probability of the structure under uncertainty. To this end, a probabilistic failure criterion is added to the constraints in the optimization problem resulting in a reliability-based design optimization (RBDO) problem \citep{enevoldsen1994reliability,kale2008tradeoff,valdebenito2010survey,lopez2012reliability}.

For most of the designed structures, the intended probability of failure is small. As a result, random sampling approaches, such as Monte Carlo simulation, to estimate the failure probability (at every optimization iteration) may lead to an impractical computational cost. Instead, the reliability index approach \citep{hasofer1974exact,madsen2006methods,melchers2018structural,haldar2000reliability} and performance measure approach \citep{tu1999new} are frequently used \citep{valdebenito2010survey}. These approaches use a first-order Taylor series approximation of the limit state function and transform the uncertain parameters to standard Gaussian random variables. However, these can introduce errors for nonlinear limit state functions and non-Gaussian random inputs. In addition, these methods typically involve two nested loops wherein the inner loop estimates the probability of failure and the outer one iterates on the design variables \citep{ramu2006inverse,acar2007reliability}. To reduce the computational cost and avoid the two-loop approach,  decoupling methods have been proposed in the past. These methods replace the inner loop of the reliability analysis with approximations \citep{yang2004experience}. For example, in \cite{du2004sequential}, the reliability constraint was replaced by a deterministic one at every iteration to sequentially update the design and estimate the reliability constraint. This deterministic constraint is then moved towards the probabilistic constraint by a shifting value obtained from a first-order approximation of the failure probability. The Karush-Kuhn-Tucker (KKT) condition can also be used \citep{kuschel1997two,kharmanda2002efficient,agarwal2007inverse} to avoid the inner reliability loop. 
\cite{cheng2006sequential} used approximations of the objective and the reliability constraint to sequentially construct optimization problems in a single-loop approach. \cite{taflanidis2008stochastic,taflanidis2008efficient} explored the design and uncertain parameter space at the same time assuming uncertainty in the design parameters. 
Further, multi-fidelity methods \citep{gano2006reliability,chaudhuri2019reusing}, response surface methods \citep{foschi2002reliability,agarwal2004reliability}, and other surrogate models \citep{missoum2007convex,zhang2004performance,bichon2008efficient,basudhar2008adaptive,suryawanshi2016reliability,moustapha2019surrogate} can be used to reduce the computational cost. 


As the estimation of failure probability remains challenging in the RBDO problems, other forms of approximations that depend on the design parameters have been used in the past for this task. 
For example, \cite{gasser1997reliability}, \cite{jensen2005design}, and \cite{jensen2007effects} used an exponential function of design parameters to approximate the probability of failure. 
\cite{ching2007approximate,ching2007local} used such local approximations and assumed the design parameters are uncertain as well. This approach uses subset simulation, an efficient sampling-based technique for low failure probabilities \citep{au2001estimation}, to obtain samples of the design and uncertain parameters from the failure region and then uses the principle of maximum entropy to estimate the parameters of the local approximation. Comparisons of robust, reliability-based, and risk-based optimization considering the cost of structural failure were performed in \cite{beck2012comparison} and \cite{beck2015comparison}. Benchmark structural design problems for RBDO were solved in \cite{aoues2010benchmark} using single-loop, double-loop, and decoupled approaches. 
Interested readers are referred to \cite{valdebenito2010survey} and \cite{lopez2012reliability} for an in-depth review of RBDO. Despite the significant progress, challenges remain when applying these approaches to large-scale problems with many design and uncertain parameters. As a result, only a few studies considering a large number of design parameters exist in RBDO literature. 



In topology optimization, the placement of material is optimized inside a design domain to optimize some performance criteria while meeting design constraints. One of the challenges of topology optimization stems from the large number of design parameters. Robust topology optimization (RTO) formulations can be used to limit the design sensitivity to the uncertainty in the material, loading, and geometry~\citep{alvarez2005minimization,guest2008structural,chen2010rtso,chen2011new,asadpoure2011robust,tootkaboni2012topology,maute2014touu,keshavarzzadeh2017topology,de2020topology}. 
Alternatively, reliability-based topology optimization (RBTO) has been performed, within a double loop strategy, to include the probability of failure as a constraint using the performance measure or reliability index approach \citep{frangopol2003life,jung2004reliability,kharmanda2004reliability,kim2006reliability,maute2003reliability,bae2002reliability,rozvany2011analytical,luo2014reliability,kang2018reliability}. Recently, \cite{da2020comparison} compared RTO to RBTO results, where reliability-based optimization is performed with a double-loop and performance measure approach. \cite{nguyen2011single} used matrix based system reliability analysis with multiple finite element mesh resolutions. \cite{silva2010component} performed RBTO in a single-loop approach using KKT condition and used both component and system failure probabilities.  \cite{jalalpour2016efficient} used second-order stochastic perturbation for estimating the sensitivity of the reliability index and used that for RBTO. 
A decoupled strategy by performing sequential optimization and reliability estimation was followed in \cite{torii2016robust} and  \cite{dos2018reliability} for stress constraints. \cite{meng2019adaptive} employed adaptively updated conjugate gradients to estimate the gradients of the probability constraints in a single-loop approach for RBTO. 
However, these works avoid the use of random sampling of the exact limit state functions and instead used the Taylor series approximation of the limit state functions, which can introduce errors in the failure probability calculations.

In this paper, we propose a stochastic gradient-based approach for RBTO. We convert the optimization problem to an unconstrained formulation that includes the reliability constraint through a penalty approach. Recently, \cite{de2020topology,de2020bi} and \cite{li2020momentum} showed that the stochastic gradient descent methods can be used to solve the RTO problems efficiently {\color{black}by using only $\mathcal{O}(1)$ random samples per optimization iteration}. Inspired by these studies, we solve the unconstrained RBTO problem using stochastic gradient descent, while focusing on reducing the computational burden of estimating the failure probability and its gradients at every iteration. In particular, we preform the estimation of the failure probability only at every few iterations using an efficient sampling strategy, such as subset simulation \citep{au2001estimation} or using surrogate models, such as polynomial expansion {\color{black}\citep{ghanem2003stochastic,doostan2011non}}. For estimating the gradients, inspired by the work of \cite{gasser1997reliability}, we use the Bayes' theorem with an exponential form for the probability density function of the design parameters in the failure region. This assumption allows for the convenient evaluation of the gradients of the failure probability. 
{\color{black}We extend this to RBTO problems, where the number of design parameters is large.
We update the parameters of this exponential form along with the design parameters using a stochastic gradient descent approach.} This is done whenever we encounter a failing design. 
As a result, this approach removes the need for the reliability index or performance measure estimates. Further, the use of stochastic gradients substantially reduces the need for a large number of realizations of the objective and constraints per optimization iteration. {\color{black}Hence, the 
proposed approach combines efficient estimation of failure probabilities with an exponential model for the probability density of design parameters from the failure region. Stochastic gradients are used to update the design parameters using only a few random samples per optimization iteration to provide computational advantage for RBTO. The proposed approach can also be extended to other RBDO problems.}

The rest of the paper is organized as follows. In Section \ref{sec:background}, we first define the topology optimization problem with a reliability constraint. We then discuss two efficient probability of failure estimation methods based on random sampling that are used in this paper. In Section \ref{sec:method}, we discuss how the stochastic gradients are used to solve the RBTO problem. We illustrate the efficacy of the proposed approach using three numerical examples thereafter. Finally, we conclude the paper with a discussion on future direction of this approach. 

\section{Background}\label{sec:background}

This section provides the setup of the RBTO problem solved in this study. Subsequently, two sampling-based methods utilized for failure probability estimation, namely subset simulation \citep{au2001estimation} and a hybrid approach with a surrogate model, will be discussed. 

\subsection{Problem Formulation}
In a deterministic optimization problem, a cost function $f(\ppm): \mathbb{R}^{\np} \rightarrow \mathbb{R}$ is minimized over the design parameters $\ppm\in \mathbb{R}^{\np}$ subject to inequality constraints $q_i(\ppm) \le 0,~~i=1,\dots,\ngg$. 
In the presence of uncertainty, the robust {\color{black}topology} optimization considers minimizing a combination of the expected value and variance of the cost function subject to a similar combination of the constraint violation \citep{beyer2007robust,de2020topology}. 
In RBTO, however, a constraint on the probability of a failure event $F$ is added to the optimization problem, e.g.,
\begin{equation} \label{eq:rbdo_defn} 
    \begin{split} 
    &\mathop{\min~}\limits_{\ppm}R(\ppm) = \Exp_{\Ym}[f(\ppm;\Ym)] \\
    &\text{subject to } C_i(\ppm) = \Exp_{\Ym}\left[q_i(\ppm;\Ym)\right] \leq 0,~~~~i=1,\dots,\ngg,\\
    & \qquad \qquad \!\! P_{F}(\ppm) = \mathbb{P}_{\xii}(F|\ppm) = \Exp_{\Ym}\left[\mathbb{I}_{F}(\Ym|\thetaa)\right] \leq p_{a},\\
\end{split} 
\end{equation} 
where $\Exp_{\xii}[\cdot]$ and $\mathrm{Var}_{\xii}(\cdot)$ denote, respectively, the expectation and variance of their arguments with respect to the probability density function (pdf) {\color{black}$p(\xii)$} of the uncertain parameters $\Ym$. 
Here, $F=\{\xii:g(\thetaa;\xii)\leq0\}$ is the failure event with limit state or performance function $g(\thetaa;\xii)$; $P_{F}(\ppm)$ is the probability of the failure event for the design parameters $\thetaa$; and $\mathbb{P}_{\xii}(\cdot)$ denotes the probability of its argument with respect to the probability measure of $\Ym$. Additionally, $\mathbb{I}_{F}(\xii|\thetaa)$ is the indicator function for the failure event $F$ for a realization of the uncertain parameters $\xii$ given a design $\thetaa$ and $p_a$ is a given maximum allowable value for $P_{F}(\ppm)$. 
\begin{figure}[!htb]
    \centering
\includegraphics[scale=1]{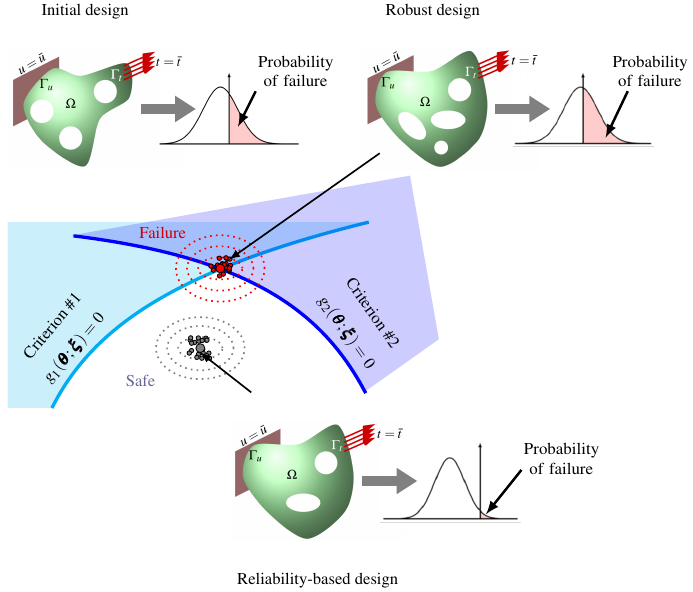}
\caption{A schematic showing that starting from an initial design, a robust design obtained by minimizing an expected value of the cost function may have a large probability of failure, when considering the failure event $F=\{\xii:g(\thetaa;\xii)=\min(g_1(\thetaa;\xii),g_2(\thetaa;\xii))\leq0\}$. 
However, the aim of the {\color{black}RBTO} is to reach a design with a prescribed (small) probability of failure. }
    \label{fig:rbdo_schem}
\end{figure} 
The design parameters may also be bounded, i.e., $\thetaa_i\in\left[\ppm_{\mathrm{min},i},  \ppm_{\mathrm{max},i}\right]$ for $i=1,\dots,\np$. 

In general, the solutions of the robust topology optimization and RBTO
lead to different designs. 
Figure \ref{fig:rbdo_schem} shows a schematic of {\color{black}RBTO} compared to the robust design approach. Here, the failure event is defined as the combination of two criteria $F=\{\xii:g(\thetaa;\xii)=\min(g_1(\thetaa;\xii),g_2(\thetaa;\xii))\leq0\}$. A robust design 
may result in a topology that has a large probability of failure. On the other hand, a reliability-based design, which is obtained by using the failure event as a constraint may have a smaller probability of failure. 



\subsection{Estimation of Failure Probability $P_{F}(\ppm)$} 

In this paper, we consider two advanced sampling methods that use conditional probability and surrogate models to efficiently estimate the failure probability $P_{F}(\ppm)$, namely, subset simulation \citep{au2001estimation} and a hybrid approach based on a surrogate model \citep{li2010evaluation}. We note, however, that other efficient sampling strategies, see, e.g., \cite{beck2015rare}, may also be employed to estimate $P_{F}(\ppm)$. Next, we briefly describe these two sampling methods.

\subsubsection{Subset Simulation} \label{sec:subsim}
Subset simulation, proposed by \cite{au2001estimation}, is a sequential Monte Carlo approach, which defines the failure region in terms of $(k+1)$ nested sets as $F:=\bigcap\limits_{j=0}^kF_j$ with $F_j:=\{\xii:g(\thetaa;\xii)\leq b_j\}$ and $b_k=0$. Hence, for these intermediary events, we have $F_j\subset F_{j-1}$ for $j=1,\dots,k$. This allows writing the failure probability as 
\begin{equation}
    P_F(\thetaa) = \mathbb{P}(F_k|\thetaa) = \mathbb{P}(F_0|\thetaa) \prod_{j=0}^{k-1} \mathbb{P}(F_{j+1}|F_{j},\thetaa).  
\end{equation} 
Since the level $b_j$ for the failure event $F_j$ is gradually reduced to zero, the first level $b_0$ is selected large enough so that $\mathbb{P}(F_0|\thetaa)$ can be easily estimated by the standard Monte Carlo sampling as 
\begin{equation}
    \mathbb{P}(F_0|\thetaa) \approx \widehat{P}_{F_0} = \frac{1}{N} \sum_{i=1}^N \mathbb{I}_{F_0} (\xii_i|\thetaa) 
\end{equation}
for $N$ independent and identically distributed (i.i.d.) samples $\{\xii_i\}_{i=1}^N$ from the pdf 
$p(\xii)$, where $N$ does not need to be very large. Next, starting with samples from $F_{j}$ the conditional probabilities are estimated as 
\begin{equation}
    \mathbb{P}(F_{j+1}|F_{j},\thetaa) \approx \widehat{P}_{F_{j+1}} = \frac{1}{N} \sum_{i=1}^N \mathbb{I}_{F_{j+1}} (\xii_i|\thetaa),\quad j=0,\dots,k-1, 
\end{equation} 
with $\{\xii_i\}_{i=1}^N$ samples from Markov chains with stationary pdf $p(\xii|F_j)$. 
Herein, we use the modified Metropolis algorithm proposed in \cite{au2001estimation} to generate these samples. 
Finally, the probability of failure is estimated as 
\begin{equation}
    \widehat{P}_F(\thetaa) = \prod_{j=0}^{k} \widehat{P}_{F_j}. 
\end{equation}
In this work, the samples are sorted by descending values of the limit state function $g(\thetaa;\xii)$, and then the levels $b_j$ are chosen as the limit state value corresponding to the $\lceil Np_0\rceil$th sample for a number $p_0$ between 0 and 1. Hence, $p_0$ can be thought of as the failure probability for each of the intermediary failure events. 
This choice of $b_j$ results in $\widehat{P}_{F_j} = p_0$ for $j=0,\dots,k-1$. For the final $k$th level and the set
\begin{equation*}
\{g(\thetaa;\xii_i): g(\thetaa;\xii_1)>g(\thetaa;\xii_2)>\dots>g(\thetaa;\xii_{N_f})\geq0>g(\thetaa;\xii_{N_f+1})>\dots>g(\thetaa;\xii_N)\}_{i=1}^N,
\end{equation*}
$\widehat{P}_{F_k} = N_f/N$, which gives the probability of failure of the current design as $\widehat{P}_F(\thetaa)=\frac{N_f}{N}p_0^{k}$. The parameter $p_0$ is generally chosen to be between 0.1 and 0.3 \citep{zuev2012bayesian}. This keeps the computational cost of estimating small failure probabilities reasonable as the number of samples required to estimate the conditional probabilities  $\mathbb{P}(F_{j+1}|F_{j},\thetaa)$ does not need to be large for a reasonable accuracy. 
Algorithm \ref{alg:SubSim} presents the steps of subset simulation. 

\begin{algorithm}[!htb]
\caption{Subset Simulation \citep{au2001estimation}}
\label{alg:SubSim}
\begin{algorithmic}
\STATE Given $N$ and $p_0$ (typically between 0.1 and 0.3)
\STATE Generate $N$ i.i.d. samples $\{\xii_i\}_{i=1}^N$ from $p(\xii)$ 
\STATE Generate $\{g(\thetaa;\xii_i)\}_{i=1}^N$ and sort its elements such that $\{g(\thetaa;\xii_i): g(\thetaa;\xii_1)>g(\thetaa;\xii_2)>\dots>g(\thetaa;\xii_N)\}_{i=1}^N$
\STATE Choose $b_0=g\left(\thetaa;\xii_{\lceil Np_0\rceil}\right)$ 
\STATE Set $j=0$
\WHILE{$b_j>0$}
\STATE Define $F_j:=\{\xii:g(\thetaa;\xii)\leq b_j\}$
\FOR {$i=1,\dots, \lceil Np_0 \rceil$}
\STATE Generate $ \lfloor 1/p_0 \rfloor$ samples from a Markov chain with stationary pdf $p(\xii|F_j)$ \\starting from the $i$th sample in $F_j$ 
\ENDFOR 
\STATE Sort these new {\color{black}$N$} samples as $\{g(\thetaa;\xii_i):g(\thetaa;\xii_1)>g(\thetaa;\xii_2)>\dots>g(\thetaa;\xii_N)\}_{i=1}^N$
\STATE Set $j=j+1$ 
\STATE Set $b_j = g\left(\thetaa;\xii_{\lceil  Np_0 \rceil }\right)$ 
\ENDWHILE 
\STATE Find the index $N_f$ such that $g\left(\thetaa;\xii_{ N_f }\right)\geq0>g\left(\thetaa;\xii_{ N_f + 1 }\right)$
\STATE Estimate probability of failure $\widehat{P}_F(\thetaa)=\frac{N_f}{N}p_0^{j}$
\end{algorithmic}
\end{algorithm}

\subsubsection{Hybrid Approach} \label{sec:hybrid}

{\color{black}The hybrid approach takes advantage of a surrogate model $\hat{g}(\thetaa;\xii)\approx g(\thetaa;\xii)$ constructed using realizations of the exact limit state function $g(\thetaa;\xii)$ \citep{li2010evaluation}. Next, }$\hat{g}(\thetaa;\xii)$ is evaluated at a sufficiently large number of samples $\{\xii_i\}_{i=1}^N$ of the uncertain parameters. The samples $\xii_i$ with $\lvert \hat{g}(\thetaa;\xii_i) \rvert \leq \gamma$, for a pre-selected tolerance parameter $\gamma$, are re-evaluated using the exact limit state function $g(\thetaa;\xii)$. The failure region is then modified to incorporate the re-evaluated limit state values as 
\begin{equation}
    F:=\{\hat{g}(\thetaa;\xii)<-\gamma\}\cup\{\{\lvert\hat{g}(\thetaa;\xii)\rvert\leq\gamma\}\cap \{g(\thetaa;\xii)<0\}\},
\end{equation} 
where the region within the tolerance limit of the surrogate model, i.e., $\{\lvert\hat{g}(\thetaa;\xii)\rvert\leq\gamma\}$ is re-evaluated using the exact limit state function and replaced with $\{g(\thetaa;\xii)<0\}$. 
Finally, the probability of failure is estimated as the fraction of the samples for which the limit state function falls in the failure region as 
\begin{equation}
    \widehat{P}_F(\thetaa)=\frac{1}{N}\sum_{i=1}^N \mathbb{I}_{F}(\xii_i|\thetaa).
\end{equation}
In the present study, we use polynomial chaos expansion (PCE) \citep{ghanem2003stochastic,xiu2002wiener} as used by  \cite{li2010evaluation} for building the surrogate model $\hat{g}(\thetaa;\xii)$. The steps of this method are illustrated in Algorithm \ref{alg:Hybrid}. The interested reader is referred to \cite{li2010evaluation} for more details about its strategy and a discussion on choosing $\gamma$. 

\begin{algorithm} [!htb]
\caption{Hybrid Approach \citep{li2010evaluation}}
\label{alg:Hybrid}
\begin{algorithmic}
\STATE Given a tolerance level $\gamma$ and sufficiently large $N$ 
\STATE Construct a surrogate model $\hat{g}(\thetaa;\xii)$ for the limit state function ${g}(\thetaa;\xii)$
\STATE Generate $N$ i.i.d. samples from $p(\xii)$ 
\FOR {$i=1,\dots,N$}
\STATE Compute $\hat{g}(\thetaa;\xii_{i})$ 
\IF {$\lvert \hat{g}(\thetaa;\xii_i) \rvert \leq \gamma$} 
\STATE Compute ${g}(\thetaa;\xii_{i})$ 
\ENDIF
\ENDFOR
\STATE Define the failure region $F:=\{\hat{g}(\thetaa;\xii)<-\gamma\}\cup\{\{\lvert\hat{g}(\thetaa;\xii)\rvert\leq\gamma\}\cap \{g(\thetaa;\xii)<0\}\}$ 
\STATE Estimated probability of failure, $\widehat{P}_F(\thetaa)=\frac{1}{N}\sum_{i=1}^N \mathbb{I}_{F}(\xii_i|\thetaa)$
\end{algorithmic}
\end{algorithm}


\section{Proposed Methodology} \label{sec:method}

Our overall strategy to solve the RBTO problem \eqref{eq:rbdo_defn} is based on the stochastic gradient descent technique. Following~\cite{de2020topology}, we consider the unconstrained formulation of \eqref{eq:rbdo_defn} with the constraints imposed via penalty parameters,
\begin{equation} \label{eq:unconst_rbdo}
    \begin{split} 
    &\mathop{\min~}\limits_{\ppm}J(\ppm) = \Exp_{\Ym}[f(\ppm;\Ym)] + \sum_{i=1}^{\ngg}\frac{\kappa_{C,i}}{2} \Exp_{\Ym}\left[({q}^+_i(\ppm;\Ym))^2\right] + \frac{\kappa_F}{2} \left[\left( \ln P_F(\thetaa) - \ln p_a \right)^+\right]^2.\\
\end{split} 
\end{equation} 
Here, $\{\kappa_{C,i}\}_{i=1}^{\ngg}$ and $\kappa_F$ are {\color{black}positive} penalty parameters used to enforce the constraints; the constraint violations are defined as ${(\cdot)}^+=0$ for $(\cdot)<0$. 
{\color{black}Note that a small value for these penalty parameters will result in a design that does not satisfy the constraints. For sufficiently large values of these penalty parameters, the solutions of \eqref{eq:unconst_rbdo} and \eqref{eq:rbdo_defn} coincide \citep{luenberger1984linear}. However, if an unnecessarily large value is selected, the convergence of the optimization process will be hindered. The numerical examples used in this paper indicate that sufficiently large values for these parameters can be used to satisfy the constraints without affecting the convergence, and a rigorous exercise to find optimum values for these parameters may not be needed. In fact, a few preliminary runs for a small number of iterations are used in those examples to select these parameters.}
The gradients of the objective $J(\thetaa)$ can then be computed using
\begin{equation}
\label{eqn:grad}
\begin{split}
    \nabla_{\thetaa} J(\thetaa) = \Exp_{\Ym}[\nabla_{\thetaa} f(\ppm;\Ym)] + \sum_{i=1}^{\ngg}{\color{black}\frac{\kappa_{C,i}}{2}} &\Exp_{\Ym}\left[\nabla_{\thetaa} (q_i^+(\ppm;\Ym))^2\right]\\ &+ \kappa_F \left( \ln P_F(\thetaa) - \ln p_a \right)^+ (\nabla_{\thetaa} \ln P_F(\thetaa)).
\end{split}
\end{equation}
The main difficulty in evaluating (\ref{eqn:grad}) is that it requires the estimation of the expected values and probability of failure using, for instance, Monte Carlo simulation, stochastic collocation {\color{black}\citep{kouri2013trust,kouri2014multilevel}, or PCE \citep{tootkaboni2012topology,keshavarzzadeh2017topology}}. The computational cost of such approaches may become prohibitive when the cost function or the constraints exhibit large variance (in the case of Monte Carlo simulation){,} or the dimension of the uncertain parameters is high (in the case of stochastic collocation or PCE). To tackle this issue, we employ a stochastic approximation of the gradients in \eqref{eqn:grad}, as discussed below. 

\subsection{Use of Stochastic Gradients} 

{\color{black}Instead of using a large number of random samples to approximate the expectations in \eqref{eqn:grad}, we generate a small sample size, e.g., $\mathcal{O}(1)$, Monte Carlo estimates of these quantities. 
The key condition here is that these estimations are independently performed from one optimization iteration to the next. 
Under certain assumptions, including strong convexity of the objective function \citep{bottou_1999}, the convergence of this approach occurs in expectation.} This approach parallels the {\it mini-batch} variants of stochastic gradient descent popularly used in training deep neural networks. At $k$th iteration, we update the design parameters as 
\begin{equation} \label{eq:sgd}
    \thetaa_{k+1} = \thetaa_k - \eta \hb_k
\end{equation} 
where $\eta$ is the step size, also known as the learning rate; and $\hb_k$ 
is a stochastic estimate of the gradient $\nabla_{\thetaa}J(\thetaa)$ {\color{black}using only $n\sim\mathcal{O}(1)$ random samples}. 
In fact, $n$ can be as small as one. 
In a previous study \citep{de2020topology}, the authors have shown that such an approach {\color{black}can efficiently solve} 
robust topology optimization problems without reliability constraints. In the presence of a reliability constraint, however, this approach can become costly as we need estimates of $P_F(\thetaa)$ and $\nabla_{\thetaa}\ln P_F(\thetaa)$. To ameliorate this potentially exorbitant computational cost, we estimate $P_F(\thetaa)$ at every $m$ iterations using subset simulation of Section \ref{sec:subsim} or the hybrid approach outlined in Section \ref{sec:hybrid}. {\color{black}Note that in \cite{de2020topology}, several variants of the standard stochastic gradient descent are studied. However, these variants use past gradient information. When applying to RBTO problems, this includes failure probability violation and its gradient for a past design. As a result, it adversely affects the convergence of the proposed RBTO approach. To avoid this, in this study, we use the {\color{black}standard stochastic gradient descent} in \eqref{eq:sgd}, where the optimization algorithm converges in our numerical examples. } 

To estimate the gradient $\nabla_{\thetaa}\ln P_F(\thetaa)$, we assume the deterministic design parameters $\thetaa$ are also uncertain with a given probability distribution (e.g., uniform in $[\thetaa_{\min},\thetaa_{\max}]$), following \cite{au2005reliability}. 
Next, we use Bayes' theorem to write
\begin{equation}
\label{eqn:bayes}
    P_F(\thetaa) = \mathbb{P}(F|\thetaa) = \frac{p(\thetaa|F)  \mathbb{P}(F)}{p(\thetaa)}. 
\end{equation}
Following \cite{gasser1997reliability} and \cite{ching2007local}, we approximate $p(\thetaa|F)$ using an exponential function as 
\begin{equation} \label{eq:pdf_approx}
    p(\thetaa|F) \approx \exp \left( -\alpha -\sum_{i=1}^{n_{\thetaa}}\beta_i \theta_i \right), 
\end{equation}
where the parameters $\alpha\in\mathbb{R}$ and $\betaa:=(\beta_1, \dots,\beta_{n_{\thetaa}})\in\mathbb{R}^{n_{\thetaa}}$ need to be estimated. {\color{black}The exponential approximation of the probability density converges to the true density in a sense of relative entropy or Kullback-Leibler distance \citep{barron1991approximation} if more polynomial terms in the design parameters are incorporated in \eqref{eq:pdf_approx} and more samples of $\thetaa$ from the failure region are used. As only local approximation of $p(\thetaa|F)$ is needed to estimate the gradients and the number of design parameters in topology optimization is very large, adding quadratic terms in \eqref{eq:pdf_approx} significantly increases the computational cost and is avoided. Here, we use more samples from the failure region as the iteration progresses to get better estimates of the parameters $\alpha$ and $\betaa$.} {\color{black}Note that other choices, such as a Gaussian mixture or generative adversarial network, can be used to approximate the probability density of design parameters from the failure region. However, a detailed investigation comparing these approximations is beyond the scope of this paper.}

For the approximation in \eqref{eq:pdf_approx} to be a pdf, it needs to satisfy the constraint 
\begin{equation} \label{eq:pdf_const}
    q_{F,0}(\alpha,\betaa) := \int_{\thetaa|F} \exp \left( -\alpha -\sum_{i=1}^{n_{\thetaa}}\beta_i \theta_i \right)\mathrm{d}\thetaa - 1 = 0. 
\end{equation}
Further, constraints to satisfy the sample mean can be added as
\begin{equation}
    q_{F,l}(\alpha,\betaa) := \int_{\thetaa|F} \theta_l \exp \left( -\alpha -\sum_{i=1}^{n_{\thetaa}}\beta_i \theta_i \right)\mathrm{d}\thetaa - \mu_{l} = 0; \quad l=1,\dots,n_{\thetaa}, 
\end{equation}
where $\mu_{l}$ is the sample mean of $l$th design parameter from the failure region. 
{\color{black}While the exponential approximation in \eqref{eq:pdf_approx} has been already used in the literature \citep{gasser1997reliability,ching2007local}, we use it for RBTO problems, where the number of design parameters is very large, and propose stochastic gradient descent to update its parameters.}
{\color{black}Note that the distribution of $\thetaa$, which is a user's choice, affects the $\ln P_F(\thetaa)$ term in \eqref{eq:unconst_rbdo} and its gradient. For a uniform distribution for $\thetaa$ between $\thetaa_{\min}$ and $\thetaa_{\max}$, $\ln p(\thetaa)$ is constant and can be ignored in the optimization process. }
Together, \eqref{eqn:bayes} and \eqref{eq:pdf_approx} lead to the simplification of the {\color{black}gradients in \eqref{eqn:grad}} as
\begin{equation}\label{eq:grad_J}
    \begin{split}
    \nabla_{\thetaa} J(\thetaa) = \Exp_{\Ym}[\nabla_{\thetaa} f(\ppm;\Ym)] + \sum_{i=1}^{\ngg}{\color{black}\frac{\kappa_{C,i}}{2}} &\Exp_{\Ym}\left[\nabla_{\thetaa} (q_i^+(\ppm;\Ym))^2\right] + \kappa_F \left( \ln P_F(\thetaa) - \ln p_a \right)^+ \betaa.
\end{split}
\end{equation}


During optimization and to satisfy the constraint in \eqref{eq:pdf_const}, whenever we encounter samples from the failure region, we solve the minimization problem
\begin{equation}
\label{eqn:albe_opt}
\mathop{\min~}\limits_{\alpha,\betaa}J_F(\alpha,\betaa) =
\frac{1}{2}\sum_{l=0}^{n_{\thetaa}}w_l q_{F,l}^2(\alpha,\betaa)
\end{equation}
using a stochastic gradient descent scheme, where $\{w_l\}_{l=0}^{n_{\thetaa}}$ are prechosen weights. {\color{black}To the best of our knowledge, this is the first study to use stochastic gradient descent for estimating these parameters and the gradients of the failure probability with respect to the design parameters. However, the use of stochastic gradient descent to solve \eqref{eqn:albe_opt} does not lead to a unique solution for $\alpha$ and $\betaa$. Further studies are needed to generate a unique exponential approximation for $p(\thetaa|F)$ while ensuring a cost that remains scalable in $\thetaa$, e.g., linear.} In particular, the gradients of  $J_F(\alpha,\betaa)$ with respect to the parameters $\alpha$ and $\betaa$ are given by 
    \begin{equation}
    \begin{split}
        &\nabla_{\!_{\alpha,\betaa}} J_{F}(\alpha,\betaa) = \sum_{l=0}^{n_{\thetaa}}w_lL_l(\alpha,\betaa) q_{F,l}(\alpha,\betaa); \quad\\
        &L_0(\alpha,\betaa) = \left[\begin{array}{c}
             -\int_{\thetaa|F} \exp \left( -\alpha -\sum_{i=1}^{n_{\thetaa}}\beta_i \theta_i \right) \mathrm{d}\thetaa \\
             -\int_{\thetaa|F} \theta_{1}\exp \left( -\alpha -\sum_{i=1}^{n_{\thetaa}}\beta_i \theta_i \right) \mathrm{d}\thetaa \\
             \vdots\\
              -\int_{\thetaa|F} \theta_{n_{\thetaa}} \exp \left( -\alpha -\sum_{i=1}^{n_{\thetaa}}\beta_i \theta_i \right) \mathrm{d}\thetaa \\
        \end{array}\right]; \quad\\
        &L_l(\alpha,\betaa) = \left[\begin{array}{c}
            -\int_{\thetaa|F} \theta_{l}\exp \left( -\alpha -\sum_{i=1}^{n_{\thetaa}}\beta_i \theta_i \right) \mathrm{d}\thetaa \\
             -\int_{\thetaa|F} \theta_{l}\theta_{1}\exp \left( -\alpha -\sum_{i=1}^{n_{\thetaa}}\beta_i \theta_i \right) \mathrm{d}\thetaa \\
             \vdots\\
              -\int_{\thetaa|F} \theta_{l}\theta_{n_{\thetaa}} \exp \left( -\alpha -\sum_{i=1}^{n_{\thetaa}}\beta_i \theta_i \right) \mathrm{d}\thetaa \\
        \end{array}\right]; \quad l=1,\dots,n_{\thetaa}.\\ 
    \end{split} 
    \end{equation} 
The stochastic -- more precisely, small sample size Monte Carlo -- approximation of $\nabla_{\!_{\alpha,\betaa}} J_{F}(\alpha,\betaa)$ is generated using {\color{black}$n_F\sim\mathcal{O}(1)$} random samples of the design parameters, $\left\{\thetaa^{(j)}\right\}_{j=1}^{n_F}$, from the failure region. Specifically, we approximate $L_l(\alpha,\betaa)$ and {\color{black}$q_{F,l}(\alpha,\betaa)$} as 
\begin{equation} \label{eq:Lq_est}
\begin{split}
    &\widehat{L}_0(\alpha,\betaa) = \left[\begin{array}{c}
             -\frac{1}{n_F}\sum_{j=1}^{n_F} \exp \left( -\alpha -\sum_{i=1}^{n_{\thetaa}}\beta_i \theta_i^{(j)} \right) \\
             -\frac{1}{n_F}\sum_{j=1}^{n_F} \theta_{1}^{(j)}\exp \left( -\alpha -\sum_{i=1}^{n_{\thetaa}}\beta_i \theta_i^{(j)} \right) \\
             \vdots\\
              -\frac{1}{n_F}\sum_{j=1}^{n_F} \theta_{n_{\thetaa}}^{(j)} \exp \left( -\alpha -\sum_{i=1}^{n_{\thetaa}}\beta_i \theta_i^{(j)} \right) \\
        \end{array}\right]; \quad\\ 
        &\widehat{L}_l(\alpha,\betaa) = \left[\begin{array}{c}
             -\frac{1}{n_F}\sum_{j=1}^{n_F} \theta_l^{(j)} \exp \left( -\alpha -\sum_{i=1}^{n_{\thetaa}}\beta_i \theta_i^{(j)} \right) \\
             -\frac{1}{n_F}\sum_{j=1}^{n_F} \theta_l^{(j)}\theta_{1}^{(j)}\exp \left( -\alpha -\sum_{i=1}^{n_{\thetaa}}\beta_i \theta_i^{(j)} \right) \\
             \vdots\\
              -\frac{1}{n_F}\sum_{j=1}^{n_F} \theta_l^{(j)}\theta_{n_{\thetaa}}^{(j)} \exp \left( -\alpha -\sum_{i=1}^{n_{\thetaa}}\beta_i \theta_i^{(j)} \right) \\
        \end{array}\right]; \quad\\ 
        &\hat{q}_{F,0}(\alpha,\betaa) = \frac{1}{n_F}\sum_{j=1}^{n_F} \exp \left( -\alpha -\sum_{i=1}^{n_{\thetaa}}\beta_i \theta_i^{(j)} \right) -1;\\
        &\hat{q}_{F,l}(\alpha,\betaa) = \frac{1}{n_F}\sum_{j=1}^{n_F} \theta_l^{(j)} \exp \left( -\alpha -\sum_{i=1}^{n_{\thetaa}}\beta_i \theta_i^{(j)} \right) -\mu_l; \quad l=1,\dots,n_{\thetaa}\\
\end{split}
\end{equation}
and update $\alpha$ and $\betaa$ via the gradient descent step 
\begin{equation}\label{eq:ab}
\begin{split}
    \left[\begin{array}{cc}
         \alpha_{k+1}  \\
         \betaa_{k+1} 
    \end{array}\right] = \left[\begin{array}{cc}
         \alpha_{k}  \\
         \betaa_{k} 
    \end{array}\right] - \eta_F \sum_{l=0}^{n_{\thetaa}}w_l\widehat{L}_l(\alpha_k,\betaa_k) \hat{q}_{F,l}(\alpha_k,\betaa_k).
\end{split}
\end{equation}
Here, $\eta_F$ is a step size parameter. As the number of samples from the failure region remains small during the initial stages of the optimization, we initially use $w_l=0$ for $l=1,\dots,n_{\thetaa}$ to avoid any convergence issue, and only use non-zero weights at the end of the optimization. Also, we choose $w_1=w_2=\dots,w_l$ to give same importance to all design parameters. 
We note that the stochastic approximations of $\nabla_{\!_{\alpha,\betaa}} J_{F}(\alpha,\betaa)$ are generated independently throughout the updates \eqref{eq:ab}. {\color{black}
Since the probability measure of $\left\{\thetaa^{(j)}\right\}_{j=1}^{n_F}$ may not exactly be the probability measure of the the design parameters given the failure event, we need to use a Radon-Nikodym derivative term $\frac{p(\thetaa|F)}{\hat{p}(\thetaa)}$ in \eqref{eq:Lq_est}, where $\hat{p}(\thetaa)$ is the probability density of the samples $\left\{\thetaa^{(j)}\right\}_{j=1}^{n_F}$. 
However, we do not write this separately and assume it is absorbed in $\eta_F$. 
} 
During optimization, one may collect the set of designs  $\left\{\thetaa^{(j)}\right\}_{j=1}^{n_F}$ from the failure region over a few iterations and then proceed to update the parameters $\alpha$ and $\betaa$ {\color{black}with $n_F>1$} . {\color{black}However, in this paper, we perform the update with the current $\thetaa$ (i.e., $n_F=1$) if it fails for any of the $n$ random samples $\{\xii_i\}_{i=1}^n$. Otherwise, we keep $\alpha$ and $\betaa$ the same.} {\color{black}Once the parameters $\betaa$ are updated, we estimate the stochastic gradients at $k$th iteration as}
\begin{equation}
    {\color{black}\hb_k = \sum_{j=1}^n \nabla_{\thetaa} f(\ppm_k;\Ym_j) + \sum_{j=1}^n \sum_{i=1}^{\ngg}\frac{\kappa_{C,i}}{2} \nabla_{\thetaa} (q_i^+(\ppm_k;\Ym_j))^2 + \kappa_F \left( \ln P_F(\thetaa_k) - \ln p_a \right)^+\betaa_{k+1}.}
\end{equation}
Algorithm \ref{alg:proposed} summarizes the steps of this proposed stochastic gradient descent method for solving the RBTO problem \eqref{eq:rbdo_defn}.  

\begin{algorithm} [!htb]
\caption{RBTO using Stochastic Gradient Descent}
\label{alg:proposed}
\begin{algorithmic}
\STATE Given step sizes $\eta$ and $\eta_F$; $m$; penalty parameters {\color{black}$\{\kappa_{C,i}\}_{i=1}^{\ngg}$} and $\kappa_F$; and $\{w_l\}_{l=0}^{n_{\thetaa}}$ 
\STATE Initial values $\thetaa_1$, $\alpha_1$, and $\betaa_1$
\FOR{$k=1,\dots$}
\IF {$k/m$ is an integer}
\STATE Use efficient sampling strategy (e.g., Algorithm \ref{alg:SubSim} or \ref{alg:Hybrid}) to estimate $\widehat{P}_F\approx\mathbb{P}(F|\thetaa)$
\ENDIF
\STATE Generate $n\sim\mathcal{O}(1)$ i.i.d. samples $\{\xii_i\}_{i=1}^n$ from $p(\xii)$ 
\STATE Estimate $\{g(\thetaa_k;\xii_i)\}_{i=1}^n$ for these samples 
\IF {$g(\thetaa_k;\xii_i) \leq 0$ for any $i\in\{1,\dots,n\}$ } 
\STATE Update the pdf parameters as $    \left[\begin{array}{cc}
         \alpha_{k+1}  \\
         \betaa_{k+1} 
    \end{array}\right] \leftarrow \left[\begin{array}{cc}
         \alpha_{k}  \\
         \betaa_{k} 
    \end{array}\right] - \eta_F \sum_{l=0}^{n_{\thetaa}}w_l\widehat{L}_l(\alpha_k,\betaa_k) \hat{q}_{F,l}(\alpha_k,\betaa_k),$ 
    \STATE $\qquad \qquad \qquad \qquad \qquad \qquad \qquad \qquad \qquad \qquad \qquad \qquad \qquad \qquad$ [see \eqref{eq:ab}] 
\ELSE
\STATE $    \left[\begin{array}{cc}
         \alpha_{k+1}  \\
         \betaa_{k+1} 
    \end{array}\right] \leftarrow \left[\begin{array}{cc}
         \alpha_{k}  \\ 
         \betaa_{k} 
    \end{array}\right]$
\ENDIF 
\STATE Estimate $\hb_k = \sum_{j=1}^n \nabla_{\thetaa} f(\ppm_k;\Ym_j) + \sum_{j=1}^n \sum_{i=1}^{\ngg}\frac{\kappa_{C,i}}{2} \nabla_{\thetaa} (q_i^+(\ppm_k;\Ym_j))^2 + \kappa_F \left( \ln \widehat{P}_F - \ln p_a \right)^+\betaa_{k+1}$ 
\STATE Update the design parameters as $\thetaa_{k+1} \leftarrow \thetaa_k - \eta \hb_k$, [see \eqref{eq:sgd}] 
\ENDFOR
\end{algorithmic}
\end{algorithm}

\subsection{Computational Cost} 

The computational cost of the proposed approach is composed of three parts. {\color{black}The first and} most computationally expensive step is to estimate the failure probability $\widehat{P}_F\approx\mathbb{P}(F|\thetaa)$, which is also the case for other reliability-based optimization methods. {\color{black}We use efficient sampling strategies, i.e., subset simulation or a hybrid approach with surrogate models, to avoid {Taylor series based approximate reliability analysis and transformation of non-Gaussian random variables}}. In addition, to further reduce the computational cost, we limit the calculation of the failure probability to every $m$ design optimization iterations. In our numerical examples, we observe setting $m$ as 25 or 50 leads to a successful design, where the optimization algorithm converges and the estimated failure probability is sufficiently accurate. {\color{black} We note that the development of fast techniques to estimate $\mathbb{P}(F|\thetaa)$, especially for small failure probabilities, is an active area of research, see, e.g., Bayesian subset simulation \citep{bect2017bayesian}, large deviation theory \citep{dematteis2019extreme,tong2020extreme,tong2020optimization}. Exploring the utility of such techniques within the proposed RBTO framework is an important future research direction.}

Second, within the employed stochastic gradient descent scheme, the estimation of the gradient {\color{black}with respect to the design parameters} is performed using $n\sim\mathcal{O}(1)$ random samples of the uncertain parameters, which drastically reduces the computational cost of gradient evaluations, as compared to methods such as the standard Monte Carlo simulation, stochastic collocation \citep{kouri2013trust,kouri2014multilevel}, or PCE \citep{tootkaboni2012topology,keshavarzzadeh2017topology}. 
{\color{black}Once the gradients are estimated, the only remaining cost is associated with updating the parameters $\thetaa$, $\alpha$, and $\betaa$ using \eqref{eq:sgd} and \eqref{eq:ab}, respectively, which is similar to other first-order optimization methods.} 
Hence, the proposed stochastic gradient descent method provides an efficient method to use random sampling for solving RBTO problems, which mostly used approximate reliability analysis in the past. 

\section{Numerical Illustrations} 

In this section, we illustrate the proposed method with three numerical examples. The first example uses a benchmark problem from \cite{rozvany2011analytical} to show the accuracy of the proposed method. Then we use a design problem of a rectangular beam and a design problem of an L-shaped beam (in two- and three-dimension), two commonly used design domain geometries in topology optimization. For these design problems, we minimize a weighted sum of compliance and mass subjected to a reliability constraint. Uncertainty is assumed in the load and material property. The results will showcase the difference between a reliability-based design and a robust design for various geometries. 

\subsection{Example I: Design of a Two-bar Truss} \label{sec:ex1}
\begin{figure}[!htb]
    \centering
    \includegraphics[scale=1]{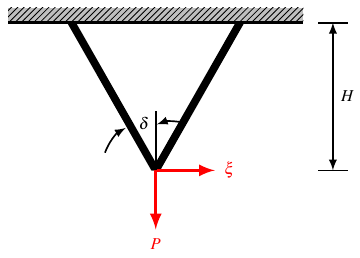}
    \caption{A benchmark problem of two-bar truss from \cite{rozvany2011analytical} used in Example I. } 
    \label{fig:benchmark}
\end{figure}

The first example uses a benchmark problem of a two-bar truss for which an analytical solution is available \citep{rozvany2011analytical}. This example is used to study the accuracy of the proposed approach as well as the influence of {\color{black}different sampling strategies to evaluate the failure probability, the penalty parameter, and the interval between two consecutive failure probability estimations} on the optimized design.  We define the problem following Section 4 of \cite{rozvany2011analytical}, where a two-bar truss is assumed with unknown cross-sectional areas and inclinations. 
Figure \ref{fig:benchmark} shows the two-bar truss with inclination $\delta\in(0,\pi/2)$ subjected to a vertical load $P$ and Gaussian distributed uncertain horizontal load $\xi$. The two bars have the same cross-sectional areas $A=\lambda A_{\max}$, where $A_{\max}$ is the maximum possible cross-sectional area and $\lambda\in[0,1]$ is a design parameter. 
The optimization problem is defined as 
\begin{equation}
\label{eqn:ex1_opt}
    \begin{split}
        &\mathop{\min}\limits_{\lambda,\delta} ~J(\lambda,\delta)=\frac{\lambda}{\cos \delta}\\
        &\text{subject to     }{\color{black}{P}_F(\lambda,\delta)} \leq p_a = 10^{-3};\\
        &\qquad \qquad ~\! 0\leq \lambda \leq 1;~~0<\delta<\pi/2,\\
    \end{split}
\end{equation} 
where the failure event $F$ is defined as $F=\{\xi:g(\lambda,\delta;\xi)=C_{\max}-C(\xi)\leq 0\}$; $C_{\max}$ is a given maximum value of the compliance; and $C(\cdot)$ is the compliance of the two-bar truss. The limit state function $g(\lambda,\delta;\xi)$ can be further simplified and written explicitly in terms of the design parameters as 
\begin{equation} \label{eq:ex1_ls}
\begin{split}
        g(\lambda,\delta;\xi) &= C_{\max} - C(\xi)\\
        &= \frac{2C_{\max}EA_{\max}}{P^2H} - \frac{1}{\lambda\cos\delta}\left( \frac{1}{\sin^2\delta} + \frac{\xi^2}{P^2\cos^2\delta} \right), \\
\end{split}
\end{equation}
where $E$ is the elastic modulus of the material and $H$ is the height of the truss as shown in Figure \ref{fig:benchmark}.

\begin{figure}[!htb]
    \centering
    \begin{subfigure}[t]{0.475\textwidth}
    \centering
    \includegraphics[scale=0.75]{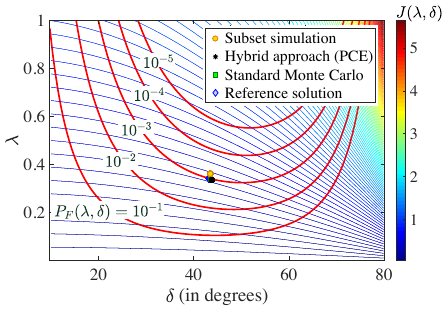}
    \caption{Comparison of different sampling strategy} \label{fig:ex1_result_a} 
    \end{subfigure} \hfill
    \begin{subfigure}[t]{0.475\textwidth}
    \centering
    \includegraphics[scale=0.75]{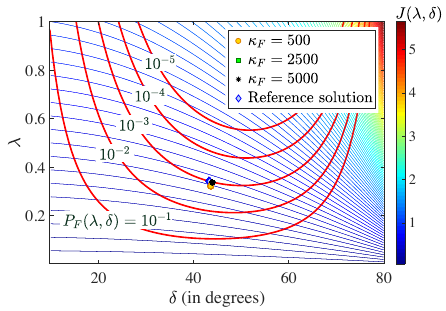}
    \caption{Results for different penalty parameter $\kappa_F$ with the hybrid approach} \label{fig:ex1_result_b} 
    \end{subfigure} \\
    \begin{subfigure}[t]{0.475\textwidth}
    \centering
    \includegraphics[scale=0.75]{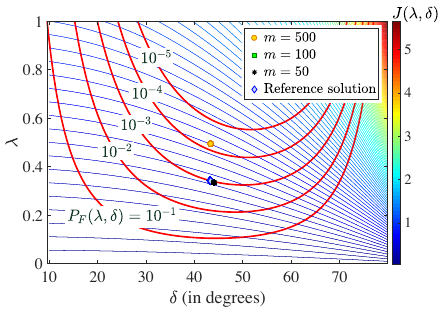}
    \caption{Results for different $m$ with the hybrid approach} \label{fig:ex1_result_c} 
    \end{subfigure}\hfill 
    \begin{subfigure}[t]{0.475\textwidth}
    \centering
    \includegraphics[scale=1.05]{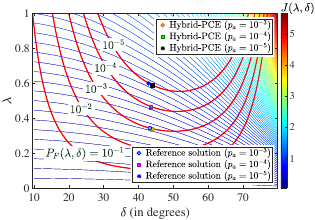}
    \caption{Results for different failure probability constraint $p_a$ with the hybrid approach} \label{fig:ex1_result_d} 
    \end{subfigure}
    \caption{Effects of different sampling methods and parameters in Example I.
    } \label{fig:ex1_results}
\end{figure}


\subsubsection{Results} 
We consider the unconstrained formulation of \eqref{eqn:ex1_opt} as in \eqref{eq:unconst_rbdo} and estimate the failure probability using standard Monte Carlo simulation or the sampling strategies discussed in Section \ref{sec:background}. We assume the maximum compliance is $C_{\max}=\frac{50P^2H}{EA_{\max}}$ and the magnitude of the vertical load $P$ is one. For the uncertain parameter $\xi$, we use a zero-mean Gaussian distribution with unit standard deviation. 
In Figure \ref{fig:ex1_results}, we plot results from different sampling strategies and parameters in the space of the optimization variables with the contours of the objective value. The five red lines are the iso-contours for the following failure probabilities $P_F(\lambda,\delta)\in\{10^{-1},10^{-2},10^{-3},10^{-4},10^{-5}\}$. The reference solution is obtained by estimating the failure probability $\mathbb{P}(F)$ using Monte Carlo simulation with $10^7$ samples, which is sufficient to estimate these failure probabilities. Table \ref{tab:ex1_results} lists the corresponding parameters and optimized design variables for all different cases considered. For all the results, we use 0.1 and $\pi/4$ for the initial values of $\lambda$ and $\delta$, respectively. The step size $\eta$ in \eqref{eq:sgd} is set to $10^{-5}$. The iterations \eqref{eq:ab} to update $\alpha$ and $\betaa$ in \eqref{eq:pdf_approx} are initialized with values $0.01$ for $\alpha$ and the entries of $\betaa$. Additionally, the step size $\eta_F$ is set to $0.2$. The reference solution is obtained after 10,000 iterations and, hence to perform a fair comparison, we run the other cases for 10,000 iterations as well. Note that the update in \eqref{eq:ab} uses $w_0=1$ and $w_l=0;l=1,\dots,n_{\thetaa}$ for the first 8000 iterations, and only for the final 2000 iterations, we use $w_0=w_1=\dots=w_{n_{\thetaa}}$ with $\sum_{i=0}^{n_{\thetaa}}w_l=1$. 

\begin{table}[!htb]
\caption{Results showing the influence of different sampling strategies and parameters of the proposed algorithm in Example I. Note that the superscript `*' denotes the optimum solution obtained for each of these cases. } \label{tab:ex1_results}
\centering 
\begin{tabular}{c c c c c c c c c} 
\hline \Tstrut 
Plot & Sampling Strategy & $\kappa_F$ & $m$ & $p_a$ & No. of $g(\lambda,\delta;\xi)$ eval. & $\lambda^*$ & $\delta^*$ & $J(\lambda^*,\delta^*)$ \\ [0.5ex] 
\hline \Tstrut 
\multirow{4}{*}{Figure \ref{fig:ex1_result_a}} & Subset simulation & 2500 & 100 & $10^{-3}$ & $1.82\times10^5$ & 0.3832 & 43.50$^\circ$ & 0.5283 \\ 
 & Hybrid approach (PCE) & 2500 & 100 & $10^{-3}$ & $7\times10^{4}$ & 0.3355 & 44.18$^\circ$ & 0.4678 \\
 & Standard Monte Carlo & 2500 & 100 & $10^{-3}$ & $10^8$ & 0.3378 & 43.79$^\circ$ & 0.4679 \\
 & Reference solution & 2500 & 1 & $10^{-3}$ & $10^{11}$ & 0.3311 & 45.00$^\circ$ & 0.4682 \\ \hline \Tstrut
 \multirow{4}{*}{Figure \ref{fig:ex1_result_b}} & Hybrid approach (PCE) & 500 & 100 & $10^{-3}$ & $7\times10^{4}$ & 0.3230 & 43.84$^\circ$ & 0.4478 \\ 
 & Hybrid approach (PCE) & 2500 & 100 & $10^{-3}$ & $7\times10^{4}$ & 0.3355 & 44.18$^\circ$ & 0.4678 \\
 & Hybrid approach (PCE) & 5000 & 100 & $10^{-3}$ & $7\times10^{4}$ & 0.3367 & 44.27$^\circ$ & 0.4702 \\
 & Reference solution & 2500 & 1 & $10^{-3}$ & $10^{11}$ & 0.3311 & 45.00$^\circ$ & 0.4682 \\ \hline \Tstrut
 \multirow{4}{*}{Figure \ref{fig:ex1_result_c}} & Hybrid approach (PCE) & 2500 & 500 & $10^{-3}$ & $2.60\times10^4$ & 0.4948 & 43.34$^\circ$ & 0.6803 \\ 
 & Hybrid approach (PCE) & 2500 & 100 & $10^{-3}$ & $7\times10^{4}$ & 0.3355 & 44.18$^\circ$ & 0.4678 \\
 & Hybrid approach (PCE) & 2500 & 50 & $10^{-3}$ & $1.24\times10^5$ & 0.3341 & 44.16$^\circ$ & 0.4657 \\
 & Reference solution & 2500 & 1 & $10^{-3}$ & $10^{11}$ & 0.3311 & 45.00$^\circ$ & 0.4682 \\ \hline \Tstrut
 \multirow{6}{*}{Figure \ref{fig:ex1_result_d}} & Hybrid approach (PCE) & 2500 & 100 & $10^{-3}$ & $7\times10^4$ & 0.3355 & 44.18$^\circ$ & 0.4678 \\ 
 & Reference solution & 2500 & 1 & $10^{-3}$ & $10^{11}$ & 0.3311 & 45.00$^\circ$ & 0.4682 \\
 & Hybrid approach (PCE) & 2500 & 100 & $10^{-4}$ & $1.14\times10^5$ & 0.4578 & 44.29$^\circ$ & 0.6396 \\
 & Reference solution & 2500 & 1 & $10^{-4}$ & $5\times10^{11}$ & 0.4651 & 43.65$^\circ$ & 0.6428 \\
 & Hybrid approach (PCE) & 2500 & 100 & $10^{-5}$ & $1.59\times10^5$ & 0.5888 & 43.99$^\circ$ & 0.8184 \\
 & Reference solution & 2500 & 1 & $10^{-5}$ & $10^{13}$ & 0.5997 & 42.88$^\circ$ & 0.8184 \\[1ex] 
\hline 
\end{tabular}
\label{table:nonlin} 
\end{table} 

Figure \ref{fig:ex1_result_a} compares the results of three sampling strategies, i.e., subset simulation, hybrid approach with PCE as the surrogate model, and standard Monte Carlo simulation, against the reference solution. In subset simulation, we use $p_0=0.1$ and, at every level, {\color{black}$N=500$} samples {\color{black}to} estimate the conditional probabilities. With the hybrid approach, during every reliability estimation, we use 100 evaluations of the limit state function to estimate the coefficients of PCE using least squares regression \citep{hadigol2018least}, and then generate $10^6$ evaluations of the PCE model to estimate the failure probabilities. The threshold $\gamma$ in Algorithm \ref{alg:Hybrid} is set to 2.5. Note that for this example building the PCE of the limit state function is computationally cheap as the dimension of the input uncertainty is one. This, however, is not the case for the following two examples. For the standard Monte Carlo simulation, during every probability of failure calculation, we use $10^6$ evaluations of $g(\lambda,\delta;\xi)$. Figure \ref{fig:ex1_result_a} shows that the hybrid approach converges to a solution very close to the reference solution and similar to the standard Monte Carlo method. However, the number of limit state function evaluations are four orders of magnitude smaller compared to using standard Monte Carlo method. With subset simulation, the optimized solution is not as accurate as the other two approaches, but the accuracy can be improved using more function evaluations to estimate the failure probability. Note that {\color{black}this} parameter setting of the subset simulation leads to about the same number {\color{black}of $g(\lambda,\delta;\xi)$ evaluation} as the hybrid approach as can be seen from Table \ref{tab:ex1_results}. Using more samples for the subset simulation approach {\color{black}would taint} the comparison among the three approaches. Table \ref{tab:ex1_results} also shows that the number of the exact limit state function evaluations for the estimation of failure probabilities is reduced by three orders of magnitude using the subset simulation or hybrid approach, compared to using a standard Monte Carlo simulation. 


For Figure \ref{fig:ex1_result_b}, we use different penalty parameters $\kappa_F$. The optimized solutions show that the accuracy improves as we use larger values. However, for larger values the convergence is slower and therefore we avoid using very large values for $\kappa_F$ in the next two examples. Similarly, Figure \ref{fig:ex1_result_c} shows that if we use large interval $m$ between two consecutive failure probability estimations (see Algorithm \ref{alg:proposed}), the accuracy of the optimized solutions deteriorates, which is expected as we are delaying the probability of failure estimation. 
Finally, we use the proposed approach for different allowable probability of failure $p_a$ ranging from $10^{-3}$ to $10^{-5}$. 
Figure \ref{fig:ex1_result_d} and Table \ref{tab:ex1_results} show that we still obtain accurate results with reasonable number of limit state function evaluations for these cases, {\color{black}when compared to the corresponding reference solutions.}

\subsection{Example II: Design of a Rectangular Beam} 

\begin{figure}
    \centering
    \begin{subfigure}[t]{\textwidth}
    \centering
    \includegraphics[scale=1]{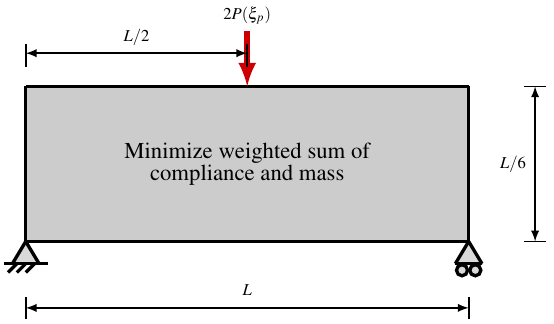}
    \caption{A rectangular beam subjected to an uncertain vertical point load at the mid span used in Example II.} 
    \label{fig:ex2_schem}
    \end{subfigure} \\
    \begin{subfigure}[t]{\textwidth}
    \centering
    \includegraphics[scale=1]{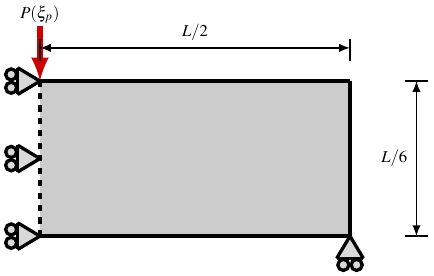}

	\caption{{\color{black}Using symmetry only right-half of the beam is considered as the design domain.}}
	\label{fig:ex2_schem_right}
	    \end{subfigure} 
	    \caption{Schematic of the rectangular beam that is designed in Example I.}
\end{figure}
In the second example, we consider the rectangular design domain shown in Figure \ref{fig:ex2_schem} subjected to an uncertain vertical load $2P({\color{black}\xi_p})$ at the mid span. We optimize the material distribution within the design domain to minimize a weighted combination of the compliance and mass subjected to a reliability constraint using the 
solid isotropic material with penalization (SIMP) approach \citep{bendsoe1989optimal,sigmund200199,sigmund2013topology}. 
A power-law model is used for interpolating the material property as follows
\begin{equation}\label{eq:simp}
E(\rho_i) = \rho_i^{\beta_\mathrm{P}} E_0;\qquad 0<\rho_i\leq 1;\qquad i=1,2,\dots,N_e,
\end{equation}
where $\rho$ is the density of a fictitious porous material; $\beta_\mathrm{P}$ is a penalization parameter; $E_0$ is the elastic modulus of the bulk material; 
and the design domain $\Omega$ is divided into $N_e$ non-overlapping elements $\{\Omega_i\}_{i=1}^{N_e}$ with volumes $\{v_i\}_{i=1}^{N_e}$. To penalize the intermediate densities, we use $\beta_\mathrm{P}=3$. The densities are defined in terms of the design parameters $\thetaa$ via a linear filter as follows 
\begin{equation}
    {\rho_i}(\thetaa) = \frac{1}{ \sum_{e=1}^{N_e}{w}_e}\sum_{e=1}^{N_e} {w}_e\theta_e; \quad i=1,\dots,N_e,
\end{equation}
where the weight ${w}_e = \max\left(0, r_f-d_{ei} \right)$ is the difference between a filter size $r_f$ and the distance $d_{ei}$ between the centers of $e$th and $i$th elements. Herein, we use 1.5 times the element width as $r_f$.
Further use of projections may be beneficial to achieve a discrete $0-1$ design \citep{sigmund2013topology}. However, we do not use any such projections in this study. 
We write the optimization problem as 
\begin{equation}\label{eq:top_obj}
\begin{split}
&\mathop{\min}\limits_{\ppm}~~J(\ppm) = \Exp\left[\sum_{i=1}^{N_e} \int_{\Omega_i} W\Big(\uu(\rho_i(\ppm);\xii),\rho_i(\ppm);\xii\Big)\mathrm{d}V_i\right] + \tau\sum_{i=1}^{N_e}v_i\rho_i(\ppm)\\
&\text{{subject to}}~~P_F(\thetaa)\leq p_a;\\ 
&\text{\phantom{subject to}}~~ {0}\leq \rho_i(\ppm)\leq {1}, \text{  for }i=1,\dots,N_e,\\
\end{split}
\end{equation} 
where the objective $J(\thetaa)$ is weighted sum of the expected value of the strain energy plus a contribution from the total mass of the structure; {\color{black}$\uu$ is the displacement vector; $\rho_i$ is the density of the $i$th element}; $W(\cdot,\cdot;\cdot)$ is the strain energy density that depends on the displacement vector $\uu$, material density $\rho$, and the uncertain variable $\Ym$; and $\tau$ is the weighting factor for contribution from the total mass to the objective. 

The uncertain load at the midspan is assumed to be $P({\color{black}\xi_p}) = P_0(1+0.25{\color{black}\xi_p})$, where ${\color{black}\xi_p}$ is a zero-mean Gaussian distributed random variable with unit standard deviation, and we assume $P_0=1$. The elastic modulus $E_0$ of the bulk material is assumed to be a log-normal random variable with unit mean and standard deviation of 0.1. 
Hence, the stochastic dimension of the problem is two. In this example, we use $\tau=0.25$ in \eqref{eq:top_obj}. 
{\color{black}We design only right-half of the beam shown in Figure \ref{fig:ex2_schem_right} using symmetry}, which we discretize into $120\times40$ bilinear elements. For the limit state function, we choose the failure as compliance value above a maximum allowable limit $C_{\max}=700$ and specify the allowable probability of failure to be $p_a=10^{-3}$. We choose these values for a scenario, where the RBTO design is different from the design that does not include the reliability constraint in the optimization process. 

\subsubsection{Results} 
\begin{figure}[!htb]
    \centering
    \begin{subfigure}[t]{0.475\textwidth}
    \centering
    {\includegraphics[scale=0.25]{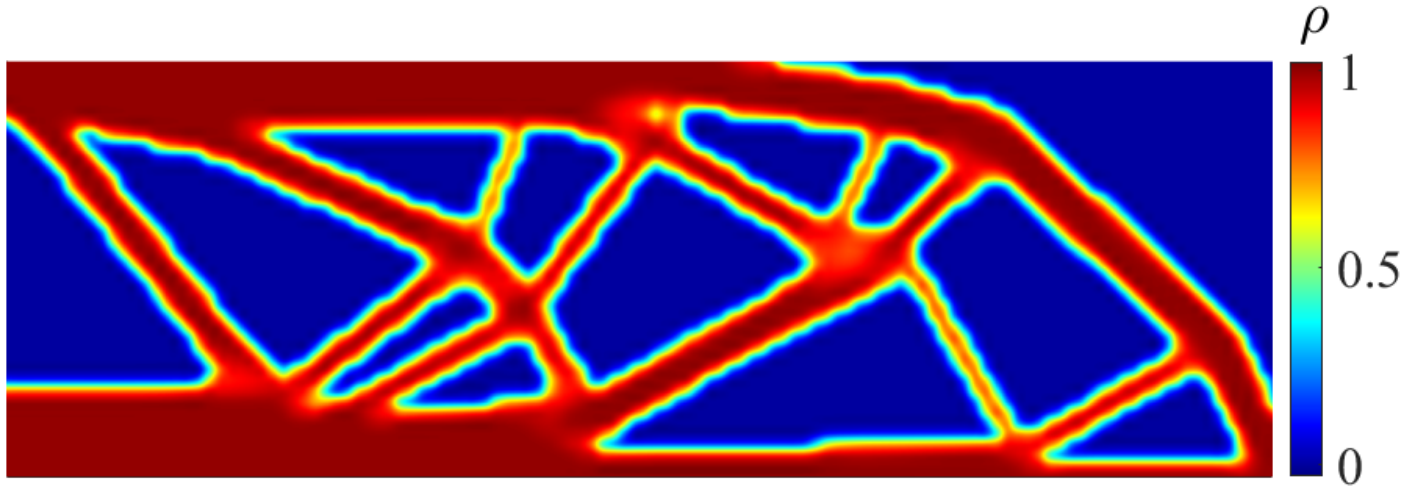}}
    \caption{Standard Monte Carlo sampling ($N=10^4$, $p_a=10^{-3}$) } \label{fig:ex2_result_a} 
    \end{subfigure} \hfill
    \begin{subfigure}[t]{0.475\textwidth}
    \centering
    {\includegraphics[scale=0.25]{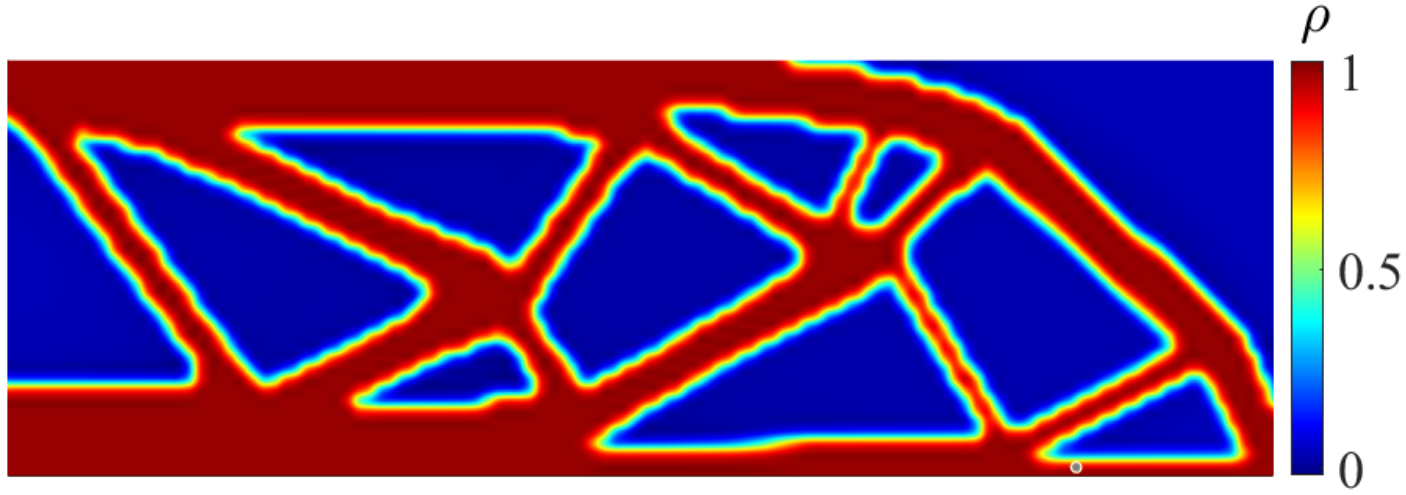}}
    \caption{Subset simulation ($p_0=0.2$, $p_a=10^{-3}$)} \label{fig:ex2_result_b} 
    \end{subfigure} 
    \begin{subfigure}[t]{0.475\textwidth}
    \centering
    {\includegraphics[scale=0.25]{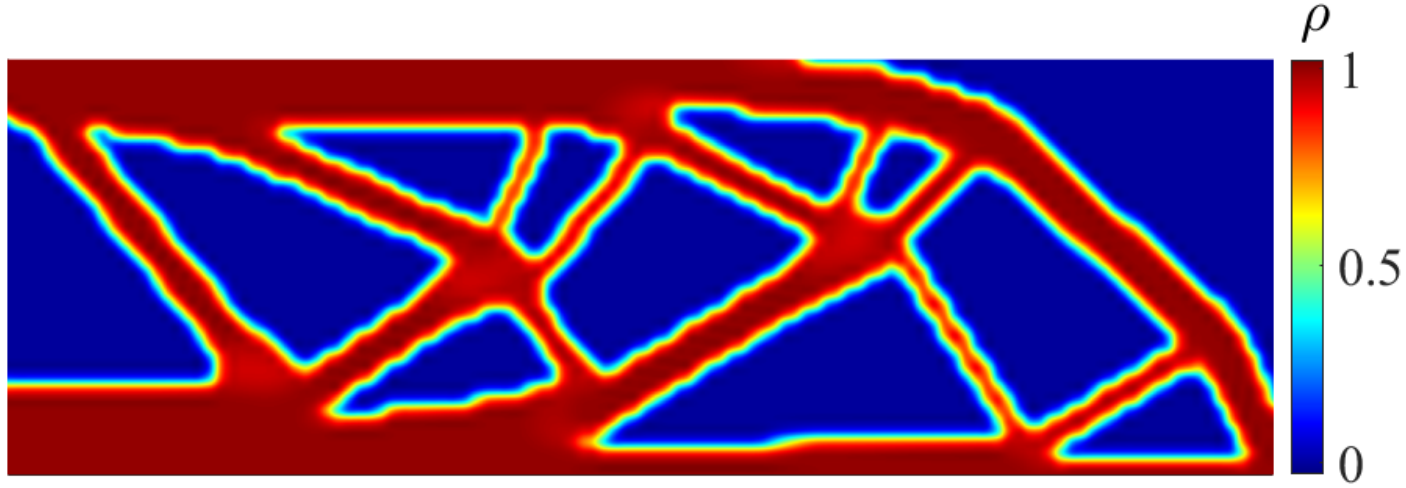}}
    \caption{Hybrid approach with PCE as surrogate model ($N=5\times10^4$, $p_a=10^{-3}$)} \label{fig:ex2_result_c} 
    \end{subfigure} \hfill 
    \begin{subfigure}[t]{0.475\textwidth}
    \centering
    {\includegraphics[scale=0.25]{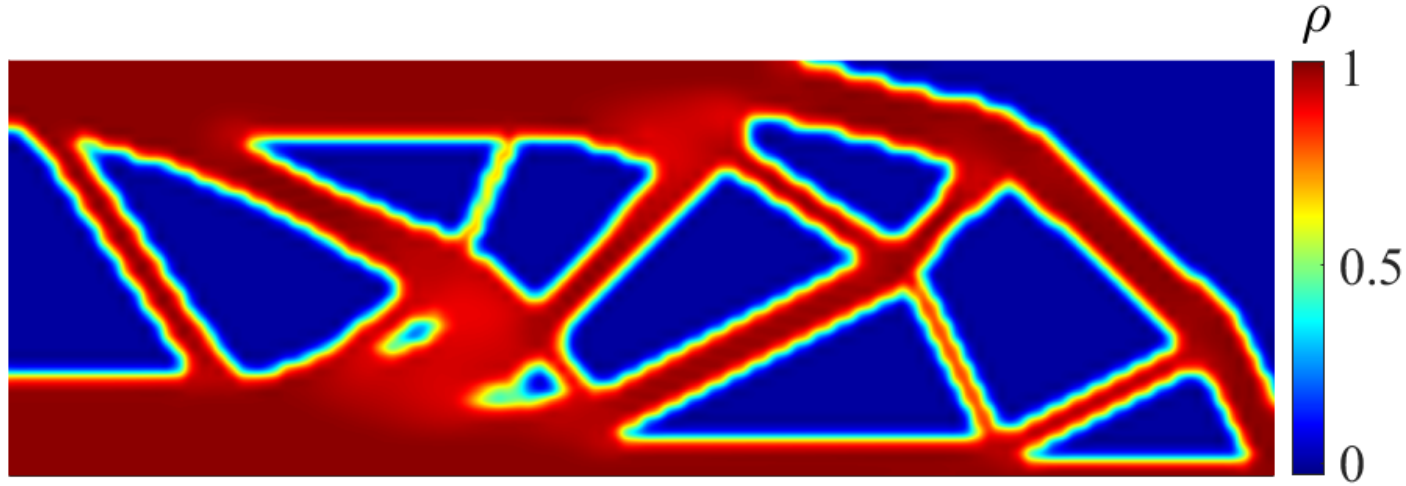}}
    \caption{Hybrid approach with PCE as surrogate model ($N=10^5$, $p_a=5\times10^{-4}$)} \label{fig:ex2_result_d} 
    \end{subfigure} 
    \begin{subfigure}[t]{0.475\textwidth}
    \centering
    {\includegraphics[scale=0.25]{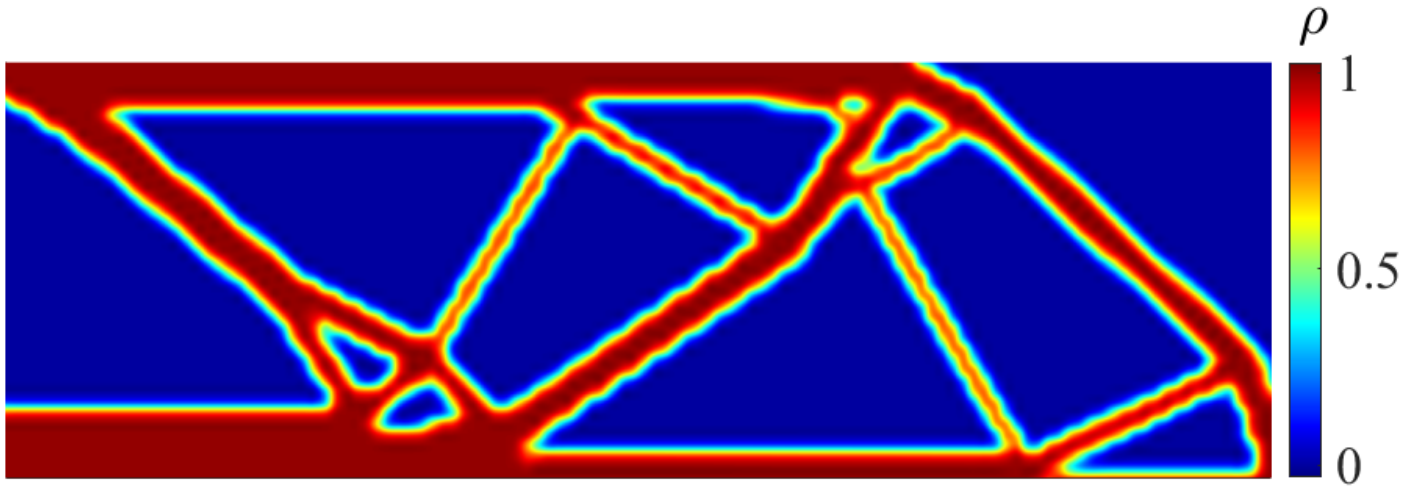}}
    \caption{Robust design that does not consider the reliability constraint} \label{fig:ex2_result_e} 
    \end{subfigure} 
    \caption{Designs obtained from different sampling strategies and for two different allowable failure probability $p_a$ in the reliability constraint in Example II. The RTO design obtained without the reliability constraint shown here for comparison. }
    \label{fig:ex2_designs}
\end{figure}

We initialize the design parameters $\thetaa$ to 0.5 each and use the stochastic gradient descent step in \eqref{eq:sgd} to perform the design optimization with $\eta=0.02$ with a mini-batch of $n=8$ random samples per iteration. We expect the parameters of the pdf $\alpha$ and $\betaa$ to be small since they need to satisfy \eqref{eq:pdf_const}. Therefore, we initialize them with $10^{-5}$ and use a small step size $\eta_F=10^{-5}$ for the updates in \eqref{eq:ab}. We perform 5000 optimization iterations, which proved sufficient for all methods to converge. Note that the update in \eqref{eq:ab} uses $w_0=1$ and $w_l=0,l=1,\dots,n_{\thetaa}$ for the first 4000 iterations, and only for the final 1000 iterations, we use $w_0=w_1=\dots=w_{n_{\thetaa}}$ with $\sum_{i=0}^{n_{\thetaa}}w_l=1$. 
Here, we use three different sampling strategies to estimate the failure probabilities at every $m=25$ iterations in the RBTO Algorithm \ref{alg:proposed}. For standard Monte Carlo sampling, we use $N=10^4$ random samples to estimate $\widehat{P}_F(\thetaa)$. For subset simulation, we set the conditional probability to $p_0=0.2$ and the number of samples to $N=1000$ for each of the levels. For the hybrid approach, we use a 4th order PCE as the surrogate model and the tolerance level of $\gamma=25$. To estimate the coefficients of the PCE surrogate, we use least squares regression and evaluate the exact limit state function for 100 realizations of the uncertain parameters. As the cost of evaluating the PCE model is negligible, we increase $N$ to $5\times10^4$. Note that these values are chosen to produce accurate estimates of the probability of failure during optimization. 
To implement the failure constraint, we set the penalty parameter to $\kappa_F=10^5$. Figures \ref{fig:ex2_result_a}, \ref{fig:ex2_result_b}, and \ref{fig:ex2_result_c} show the designs obtained from these three sampling strategies. These designs look similar except for a few extra members near left side in Monte Carlo and hybrid approach. The subset simulation design has thicker members in those places. 
The total number of finite element solves required by each of these three sampling strategies during the optimization, however, varies by more than one order of magnitude. For example, the standard Monte Carlo approach requires $\sim 2.04\times10^6$  finite element solves compared to $\sim 6.0\times10^5$ for subset simulation and $\sim6.9\times10^4$ finite element solves for the hybrid approach.



\begin{table}[!htb] 
\caption{Probability of failure $\widehat{P}_F(\thetaa)$ estimated using $10^5$ random samples and {\color{black}mass ratio} of the final designed structures shown in Figure \ref{fig:ex2_designs} in Example II. 
} 
\centering 
\begin{tabular}{c c c c c} 
\hline 
\Tstrut Design & Sampling strategy & $p_a$ & $\widehat{P}_F(\thetaa)$ & {\color{black}Mass ratio} \\ [0.5ex] 
\hline 
\Tstrut \multirow{4}{*}{Reliability-based} & Standard Monte Carlo & $10^{-3}$ & $1.7\times10^{-3}$ & 0.4617 \\ 
 & Subset simulation & $10^{-3}$ & $1.2\times10^{-3}$ & 0.4705 \\
 & Hybrid approach (PCE) & $10^{-3}$ & $1.2\times10^{-3}$ & 0.4665 \\
 & Hybrid approach (PCE) & $5\times10^{-4}$ & $5.7\times10^{-4}$ & 0.5062\Bstrut\\\hline\Tstrut
Robust & -- & -- & $2.12\times10^{-2}$ & 0.3444 \\ [1ex] 
\hline 
\end{tabular}
\label{table:ex2_Pf} 
\end{table}

The failure probabilities of these designed structures are then estimated by Monte Carlo sampling using $10^5$ evaluations of the limit state function $g(\thetaa;\xii)$. Table \ref{table:ex2_Pf} shows that the final designs obtained using these three sampling strategies have failure probabilities slightly over $10^{-3}$. This is due to the penalty formulation used here. To reach a design with probability of failure strictly below or equal to $p_a=10^{-3}$, we would have to increase the penalty parameter $\kappa_F$. However, as the results for the problem in Section \ref{sec:ex1} have shown this will likely slow down the convergence and require a large number of iterations to converge to a $0-1$ design. Instead, we can also use a smaller $p_a$ in the optimization problem than required by the engineering application. For example, we use $p_a=5\times10^{-4}$ and implement the hybrid approach with $N=10^5$ evaluations of a 4th order PCE as the surrogate model and same tolerance level $\gamma=25$ as before. The resulting design is shown in Figure \ref{fig:ex2_result_d}, which has more members as expected. The probability of failure of this designed structure estimated from $10^5$ evaluations of the limit state function $g(\thetaa;\xii)$ is $5.7\times10^{-4}$ (see Table \ref{table:ex2_Pf}). {\color{black}Table \ref{table:ex2_Pf} also lists the mass ratio for these designs, where the mass ratio is defined as the ratio of the mass of the designed structure to the mass of a structure that occupies the entire design domain. }

\begin{figure}[!htb]
    \centering
    \begin{subfigure}[t]{\textwidth}
    \centering
    {\includegraphics[scale=0.275]{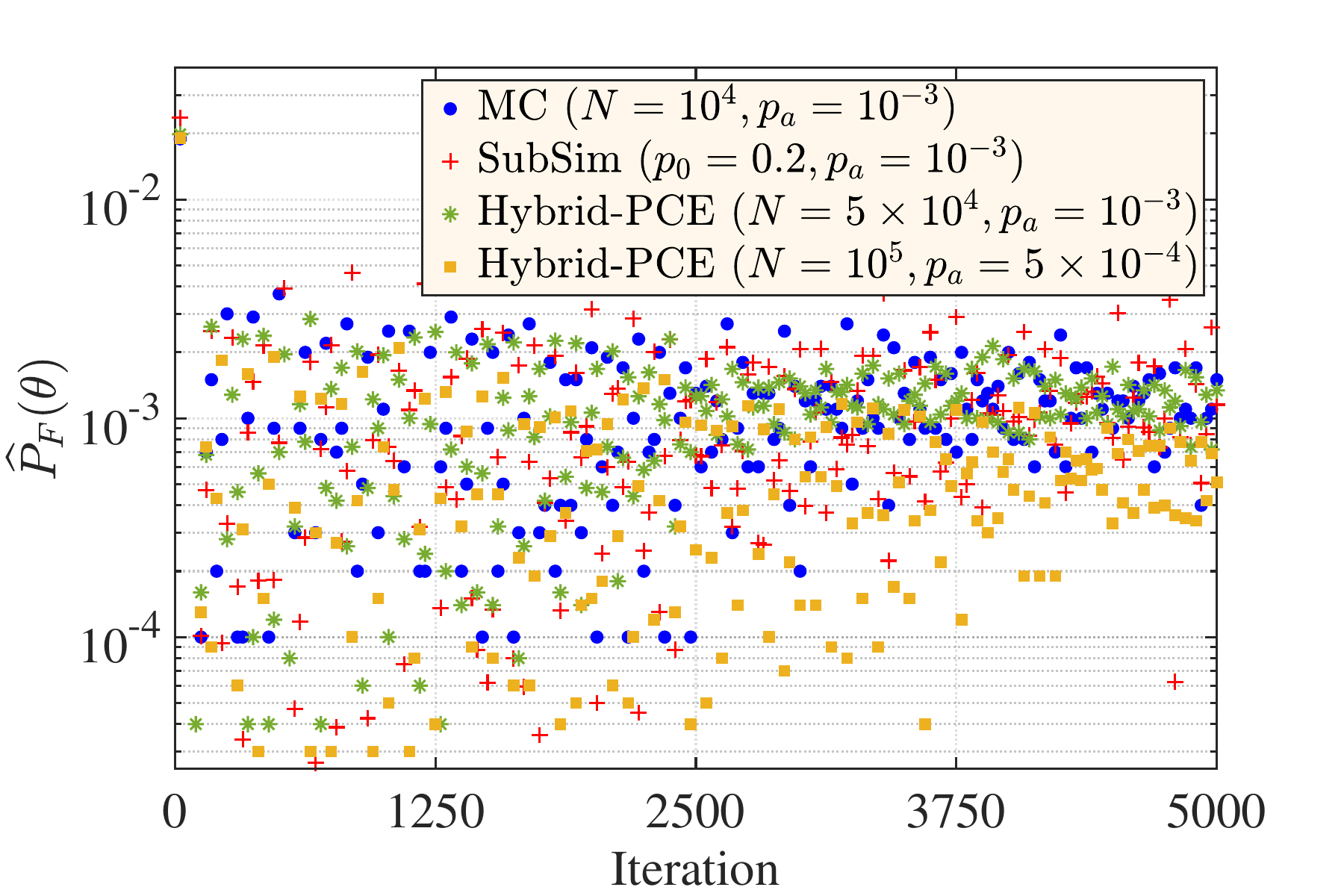}}
    \caption{ The failure probabilities of the design $\widehat{P}_F(\thetaa)$ as estimated during the optimization} \label{fig:ex2_PfHist} 
    \end{subfigure} \\
    \begin{subfigure}[t]{\textwidth}
    \centering
    {\includegraphics[scale=0.275]{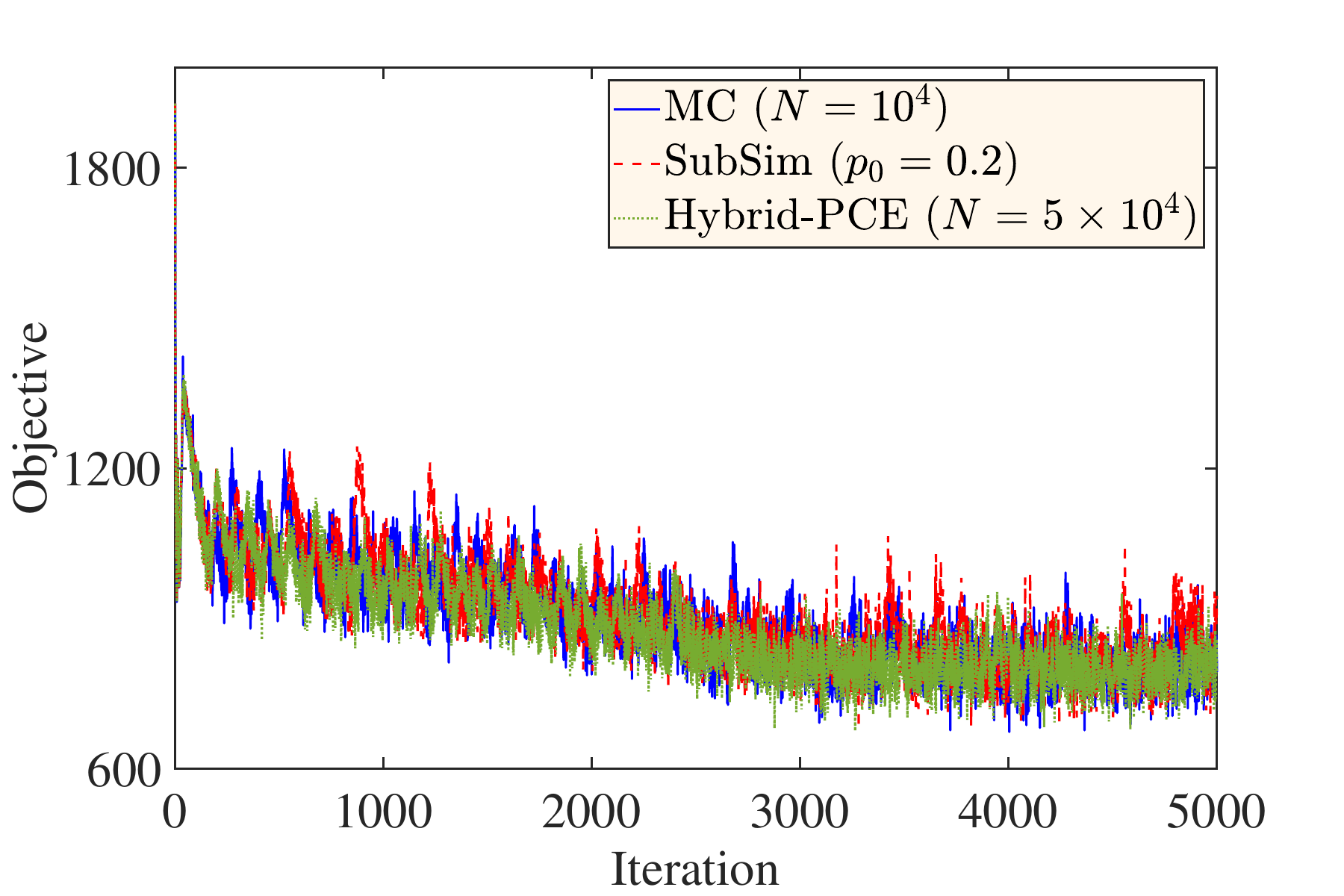}}
    \caption{Objective values for $p_a=10^{-3}$ in the reliability constraint} \label{fig:ex2_ObjHist} 
    \end{subfigure} 
    \caption{Failure probabilities and objective values estimated during the optimization for different sampling strategies in Example II. Here, MC stands for Monte Carlo method; SubSim stands for subset simulation; and Hybrid-PCE stands for the hybrid approach with PCE as the surrogate model.}
    \label{fig:ex2_results}
\end{figure} 

Figure \ref{fig:ex2_PfHist} plots the estimates of failure probability at every 25th iteration during the optimization process. This figure shows that the probability of failure initially oscillates but it starts to converge beyond 3750 iterations. Also, the plot of objective values in Figure \ref{fig:ex2_ObjHist} attests to this. 

Figure \ref{fig:ex2_result_e} shows a design obtained by solving the same optimization problem but ignoring the reliability constraint. 
The reliability-based designs are significantly different compared to this robust design. Also, the robust design has a failure probability of $2.12\times10^{-2}$, significantly higher than the allowable limit. {\color{black}Note that the objective has a contribution from the mass of the structure. As a result, the robust design has slender members compared to the reliability-based designs in the absence of a constraint on the probability of failure. This is also evident from the mass ratio of these designs reported in Table \ref{table:ex2_Pf}. }
Hence, this example shows the effectiveness of the proposed approach to produce a design that is reliable.


\subsection{Example III: Design of an L-shaped Beam} 

\begin{figure}[!htb]
    \centering


    
    
    
    

\includegraphics[scale=1]{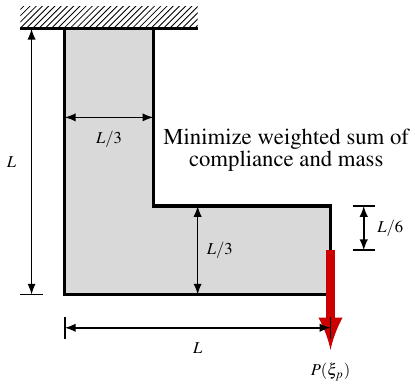}
    \caption{A two-dimensional L-shaped beam with uncertain material property and subjected to an uncertain load $P({\color{black}\xi_p})$ at the right side is used in Example III (a). } 
    \label{fig:ex3_schem}
\end{figure}

In the third example, we seek to find a structure within an L-shaped design domain subjected to material and loading uncertainty. Figures \ref{fig:ex3_schem} and \ref{fig:ex3_schem3D} show schematics of the problem in two and three dimensions, respectively, where an uncertain vertical load $P({\color{black}\xi_p})=P_0(1+0.5{\color{black}\xi_p})$ is applied at the center of the most right vertical edge/face of the design domain. We model ${\color{black}\xi_p}$ as a zero-mean Gaussian random variable with unit standard deviation and set $P_0=0.5$. The elastic modulus $E_0$ of the bulk material is assumed to be {\color{black}a lognormal random variable with unit mean and standard deviation of 0.2}. 
Hence, we have two random variables with known probability distribution functions. This example is used to confirm the findings from the previous example as well as to extend the proposed approach to a three dimensional design problem. 
The optimization results for the two cases are discussed next. 

\subsubsection{Case (a): Two-dimensional Beam} 

\begin{figure}[!htb]
    \centering
    \begin{subfigure}[t]{0.475\textwidth}
    \centering
    {\includegraphics[scale=0.225]{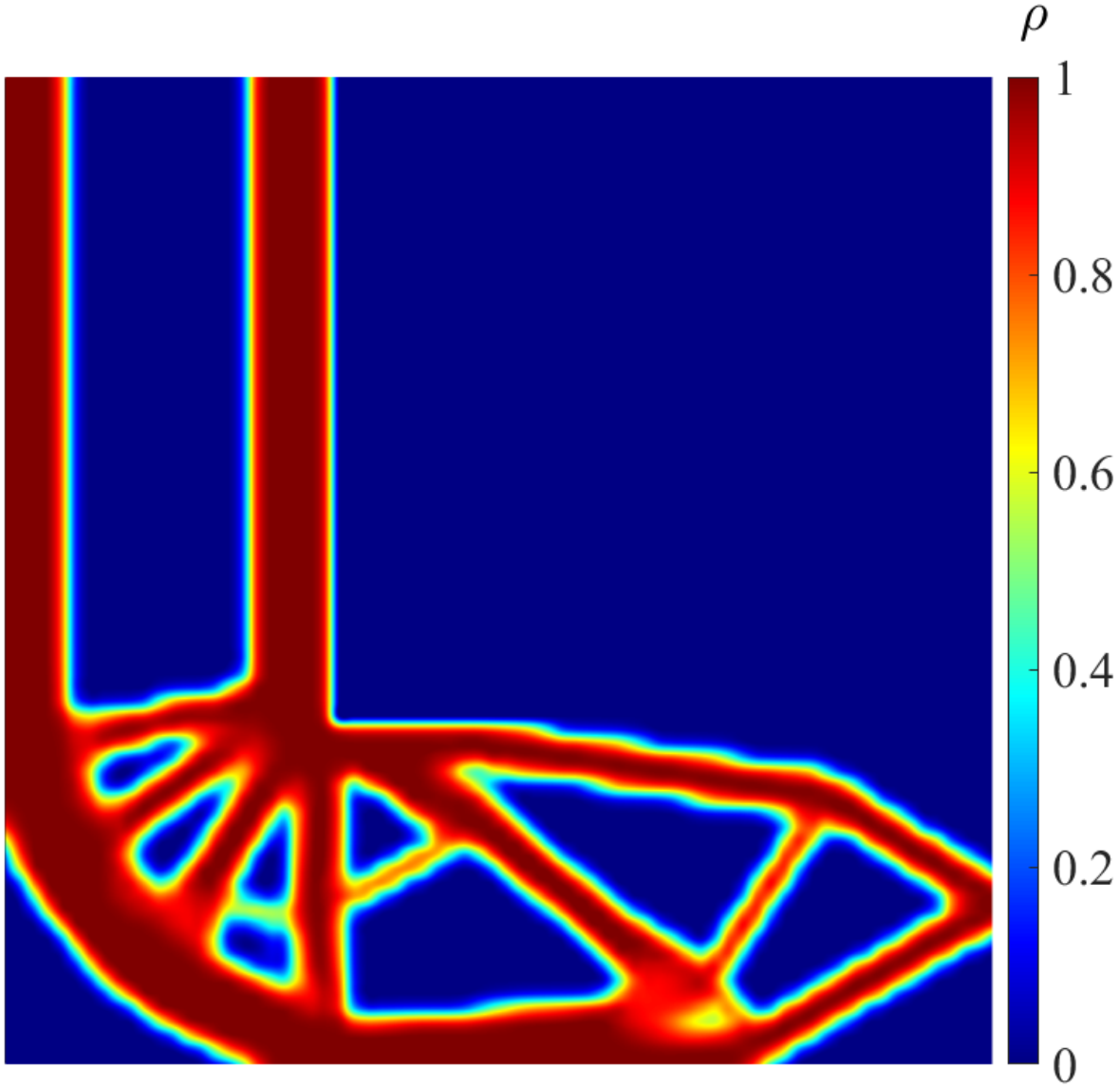}}
    \caption{Standard Monte Carlo sampling ($N=10^4$, $p_a=10^{-3}$) } \label{fig:ex3_result_a} 
    \end{subfigure} \hfill
    \begin{subfigure}[t]{0.475\textwidth}
    \centering
    {\includegraphics[scale=0.225]{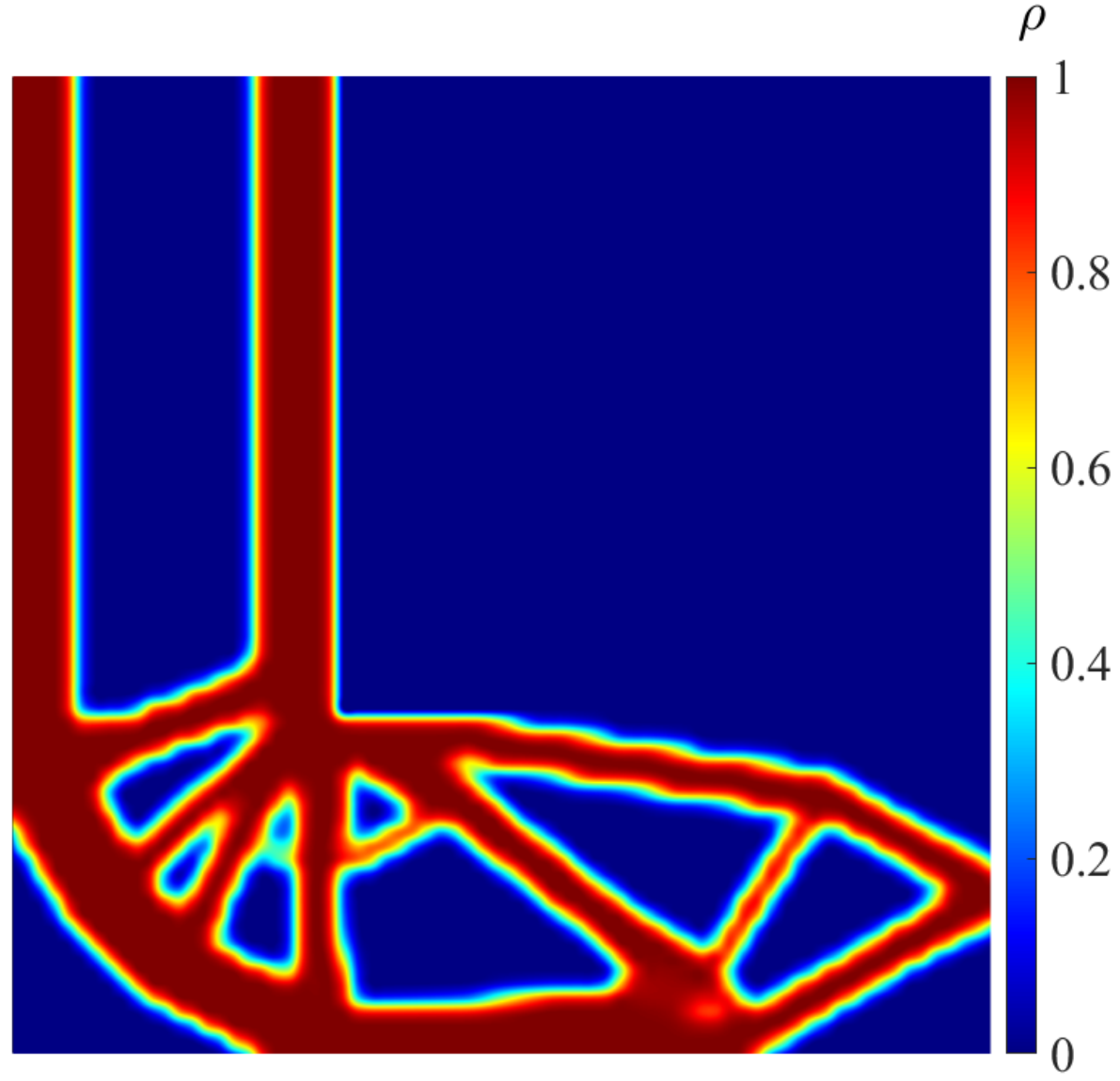}}
    \caption{Subset simulation ($p_0=0.2$, $p_a=10^{-3}$) \phantom{extra text}
    } \label{fig:ex3_result_b} 
    \end{subfigure} 
    \begin{subfigure}[t]{0.475\textwidth}
    \centering
    {\includegraphics[scale=0.225]{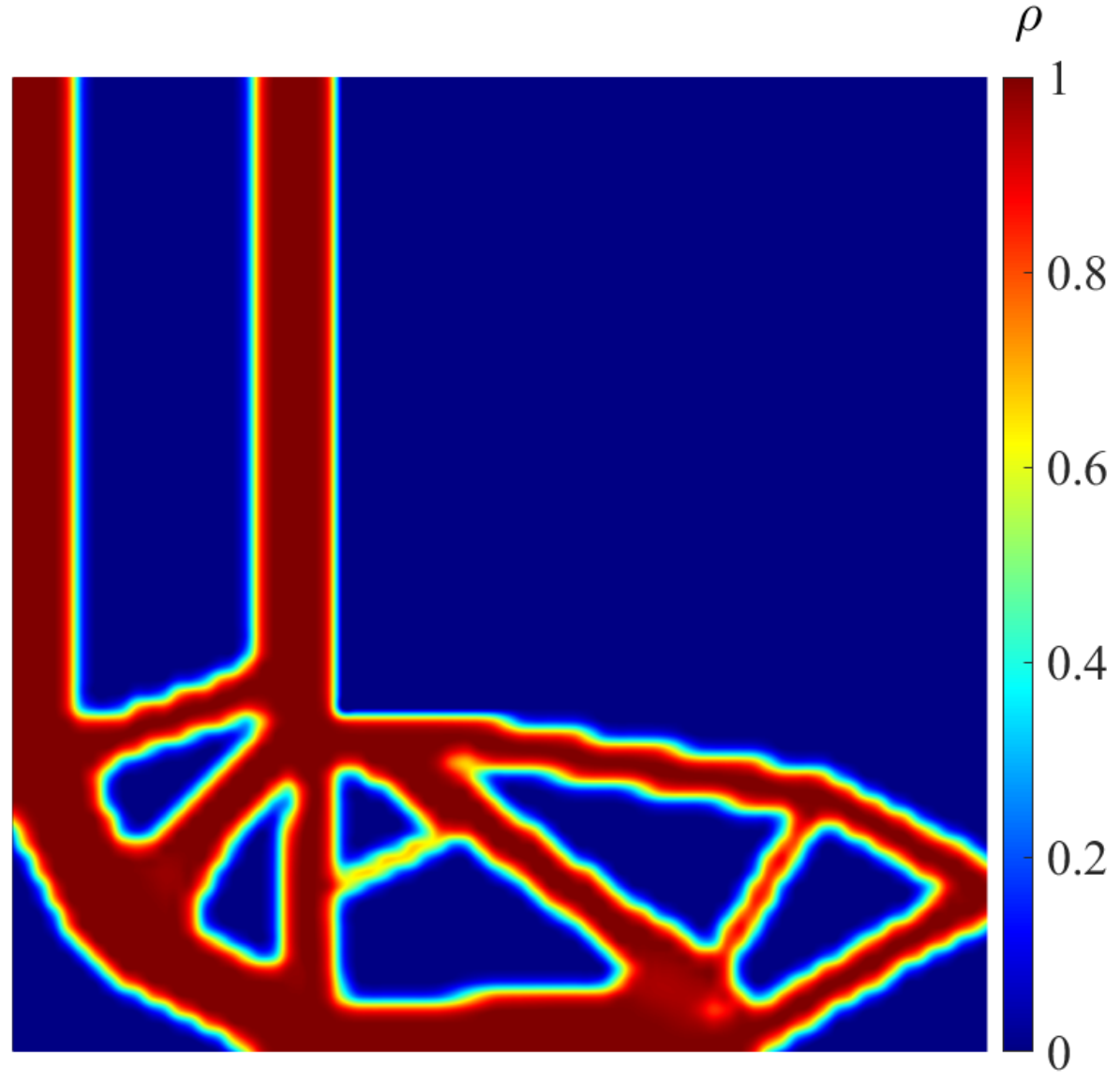}}
    \caption{Hybrid approach with PCE as surrogate model ($N=5\times10^4$, $p_a=10^{-3}$)} \label{fig:ex3_result_c} 
    \end{subfigure} \hfill 
    \begin{subfigure}[t]{0.475\textwidth}
    \centering
    {\includegraphics[scale=0.225]{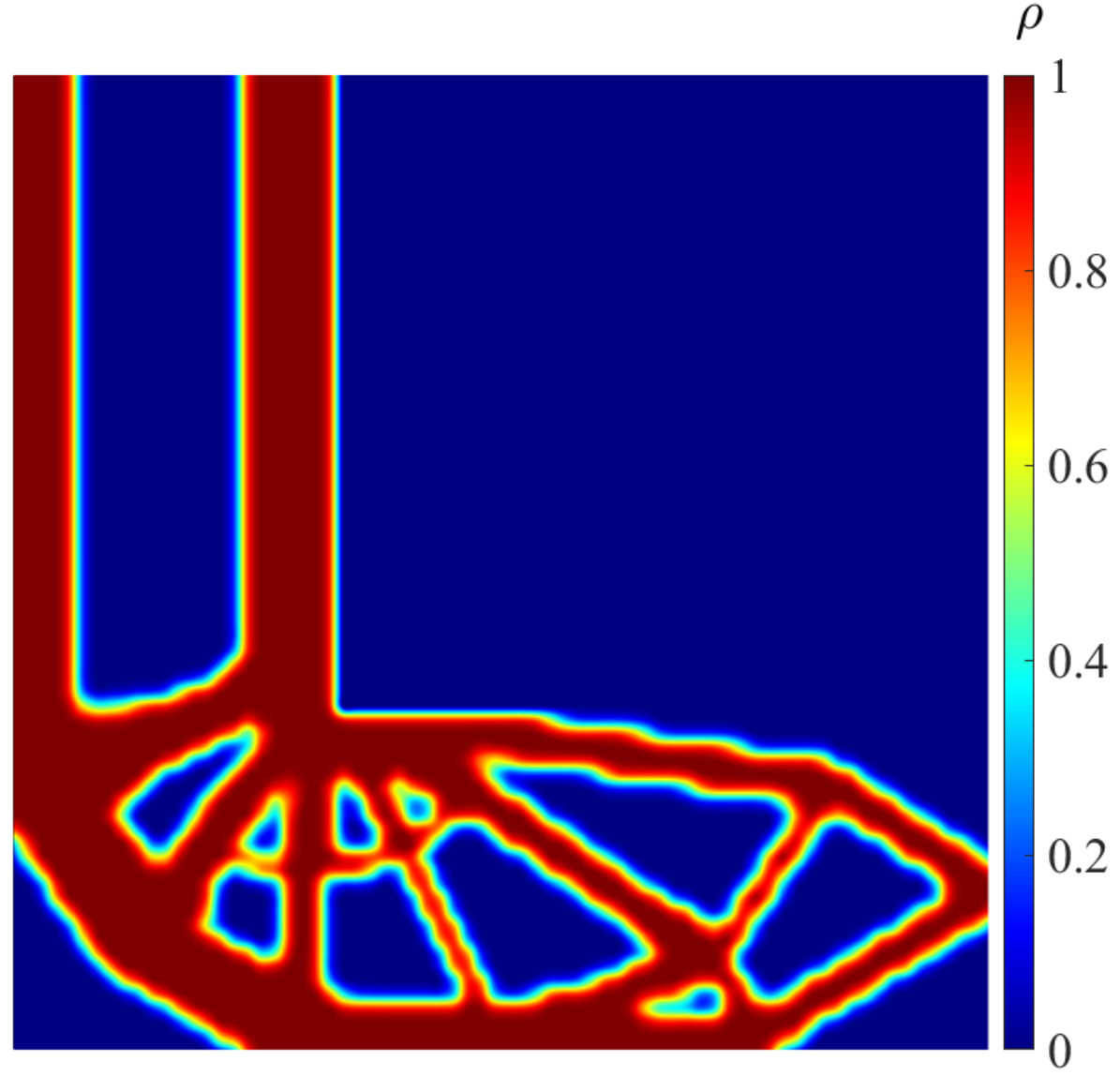}}
    \caption{Hybrid approach with PCE as surrogate model ($N=10^5$, $p_a=5\times10^{-4}$)} \label{fig:ex3_result_d} 
    \end{subfigure} 
    \begin{subfigure}[t]{0.475\textwidth}
    \centering
    {\includegraphics[scale=0.225]{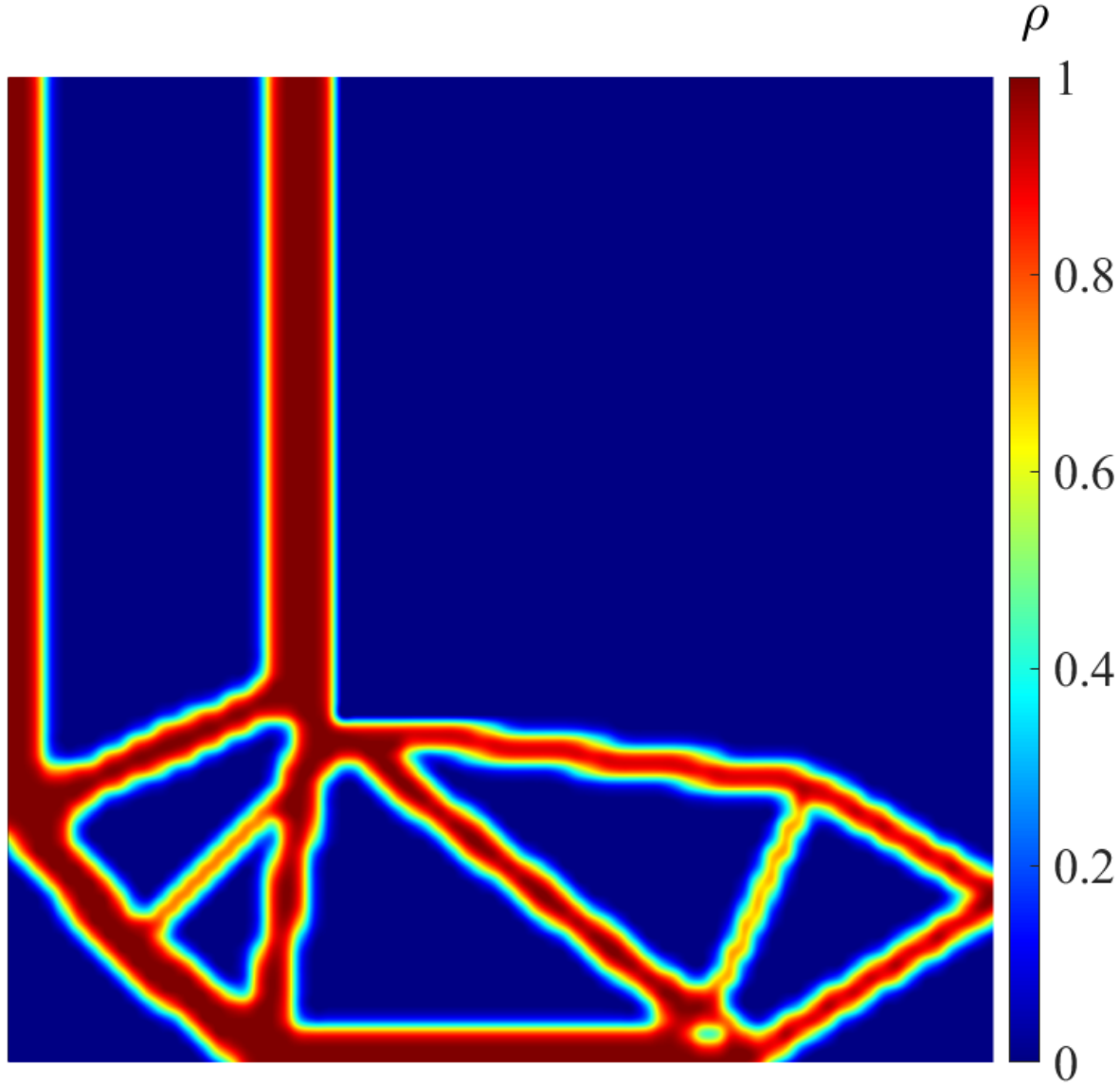}}
    \caption{Robust design that does not consider the reliability constraint } \label{fig:ex3_result_e} 
    \end{subfigure} 
    \caption{Designs obtained from different sampling strategies and for two different allowable failure probability $p_a$ in the reliability constraint in Example III (a). The RTO design obtained without the reliability constraint shown here for comparison. }
    \label{fig:ex3_designs}
\end{figure}

We solve the optimization problem as defined in \eqref{eq:top_obj} with $\tau=0.25$. We discretize the design domain into a total of 2880 bilinear elements. 
We define the failure as the compliance exceeding $C_{\max}=650$, and the reliability constraint is defined to keep the probability of failure below $p_a=10^{-3}$. 
We solve the optimization problem with a step size $\eta=0.035$ for updating the design parameters $\thetaa$ for $5000$ iterations and a mini-batch of $n=4$ random samples per iteration for the stochastic gradient step in \eqref{eq:sgd}. 
The optimization variables are initialized as in Example II. 
For the parameters $\alpha$ and $\betaa$ of $p(\thetaa|F)$ we use the same initialization and step size as in the previous example. 
Note that the update in \eqref{eq:ab} uses $w_0=1$ and $w_l=0,l=1,\dots,n_{\thetaa}$ for the first 4000 iterations, and only for the final 1000 iterations, we use $w_0=w_1=\dots=w_{n_{\thetaa}}$ with $\sum_{i=0}^{n_{\thetaa}}w_l=1$. 
Three different sampling strategies, namely standard Monte Carlo, subset simulation, and the hybrid approach with PCE as surrogate model, are used to estimate the failure probabilities at every $m=25$ iterations. To enforce the reliability constraint, we use a penalty parameter $\kappa_F=5\times10^4$. The settings for these three sampling strategies are the same as those used in Example II except for the subset simulation, where we increase $N$ to 2000 to estimate the conditional probabilities at each level. Note that these settings for the sampling strategies are selected to produce accurate estimates of the failure probabilities during optimization. Figures \ref{fig:ex3_result_a}, \ref{fig:ex3_result_b}, and \ref{fig:ex3_result_c} show the designs obtained from these three sampling strategies. The designs obtained by the standard Monte Carlo and subset simulation look similar. The design from the hybrid approach differs slightly from the other two. 
For the Monte Carlo approach the total number of finite element solves required during the entire optimization process is $\sim 2.04\times10^6$  compared to $\sim 1.60\times10^6$ for the subset simulation. However, in the hybrid approach we only perform $\sim4.8\times10^4$ solves. Hence, the hybrid approach proves to be the most economic option here. We use this approach for the three dimensional extension of this problem in the next section.

\begin{figure}[!htb]
    \centering
    \begin{subfigure}[t]{\textwidth}
    \centering
    {\includegraphics[scale=0.3]{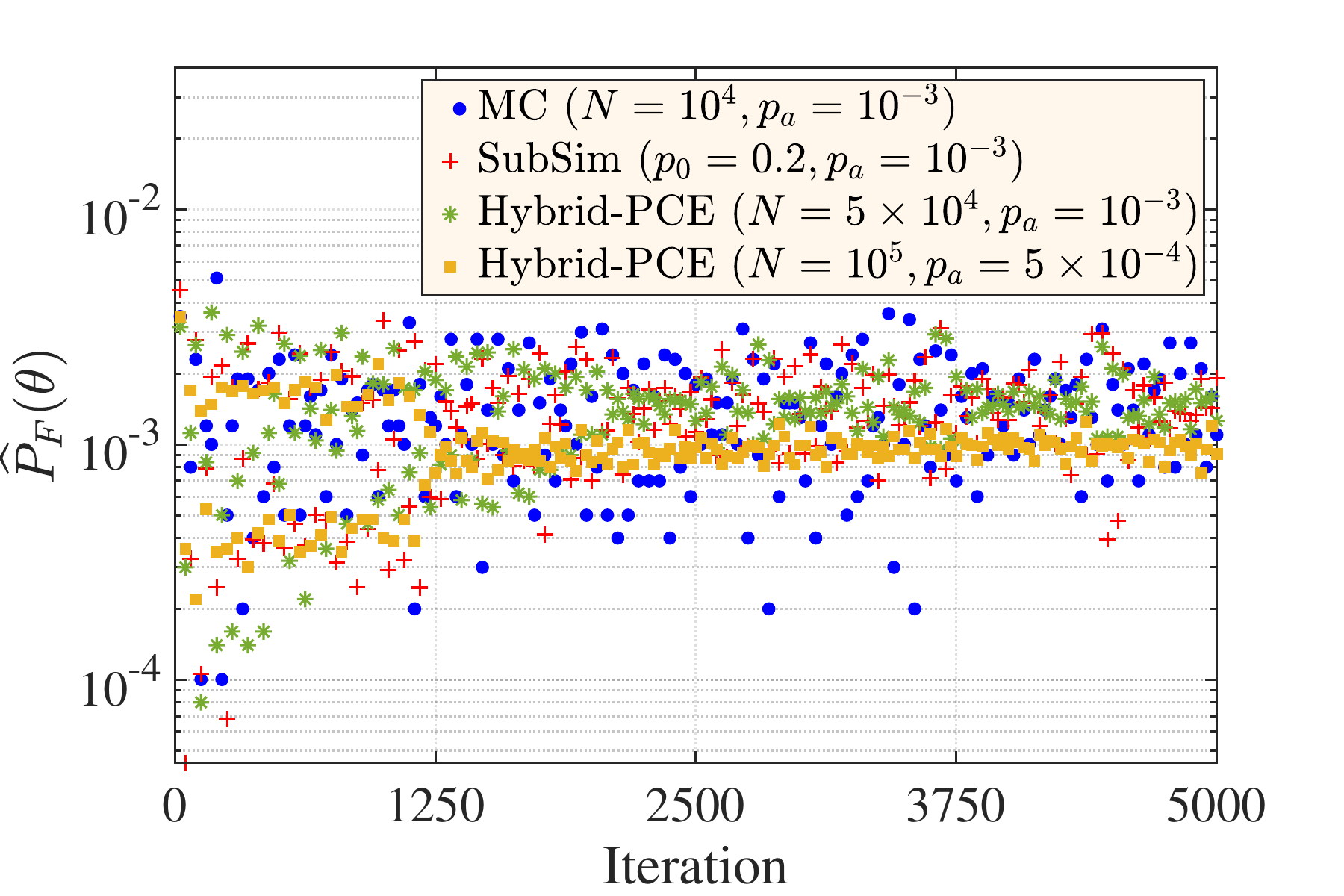}}
    \caption{ The failure probabilities of the design $\widehat{P}_F(\thetaa)$ as estimated during the optimization process} \label{fig:ex3_PfHist} 
    \end{subfigure} \\
    \begin{subfigure}[t]{\textwidth}
    \centering
    {\includegraphics[scale=0.3]{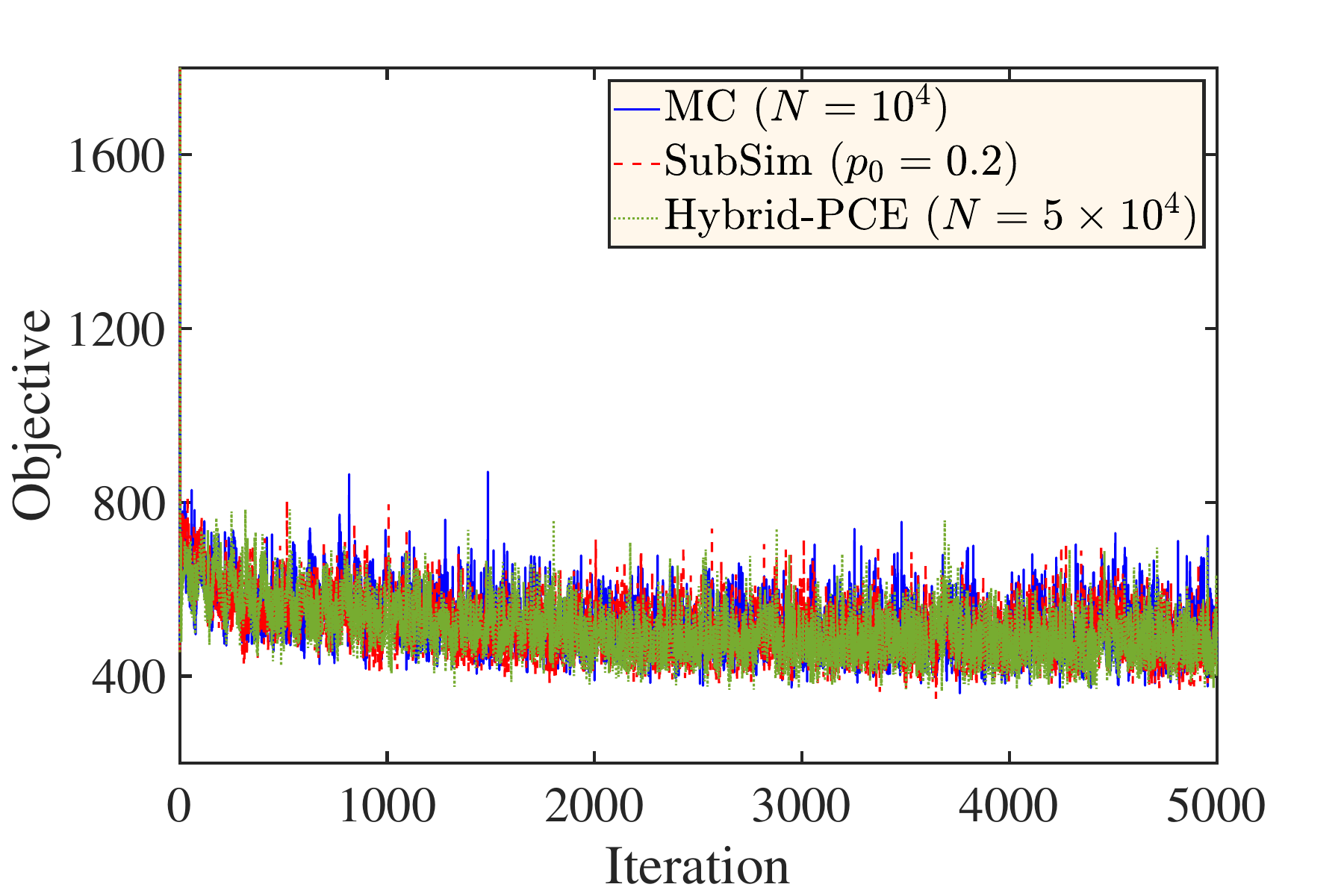}}
    \caption{Objective values for $p_a=10^{-3}$ in the reliability constraint} \label{fig:ex3_ObjHist} 
    \end{subfigure} 
    \caption{Failure probabilities and objective values estimated during the optimization process for different sampling strategies in Example III (a). Here, MC stands for Monte Carlo method; SubSim stands for subset simulation; and Hybrid-PCE stands for the hybrid approach with PCE as the surrogate model.} 
    \label{fig:ex3_results}
\end{figure} 

Table \ref{table:ex3_Pf} shows the probability of failure for these designs estimated using $10^5$ random samples. 
Designs from the three sampling strategies have failure probabilities slightly more than the prescribed value of $10^{-3}$. 
However, similar to the previous example, using the hybrid approach with a smaller allowable failure probability as $p_a=5\times10^{-4}$ and $N=10^5,\kappa_F=3\times10^4$, we obtain a design shown in Figure \ref{fig:ex3_result_d} that has more bars and a probability of failure $9.6\times10^{-4}$ as shown in Table \ref{table:ex3_Pf}. 
{\color{black}Note that the mass ratio in Figure \ref{fig:ex3_result_d} for the hybrid approach with $p_a=5\times10^{-4}$ and $N=10^5$ is largest, where the mass ratio is defined as before. We also compare the results with a RTO design that does not include the reliability constraint, but the objective formulation remains the same as in the RBTO problem. For this RTO design, the mass ratio is the lowest. 
}
The convergence of the objective and failure probabilities of the design are shown in Figure \ref{fig:ex3_results}. 
This example again shows the efficacy of the proposed approach for RBTO. In the next case, we extend the RBTO problem to the design of a three-dimensional L-shaped beam.

\begin{table}[!htb]
\caption{Probability of failure $\widehat{P}_F(\thetaa)$ estimated using $10^5$ random samples and {\color{black}mass ratio} of the final designed structures shown in Figure \ref{fig:ex2_designs} in Example III (a). 
} 
\centering 
\begin{tabular}{c c c c c} 
\hline 
\Tstrut Design & Sampling strategy & $p_a$ & $\widehat{P}_F(\thetaa)$ & {\color{red}Mass ratio} \\ [0.5ex] 
\hline 
\Tstrut \multirow{4}{*}{Reliability-based} & Standard Monte Carlo & $10^{-3}$ & $1.7\times10^{-3}$ & 0.4847 \\ 
 & Subset simulation & $10^{-3}$ & $1.4\times10^{-3}$ & 0.4971 \\
 & Hybrid approach (PCE) & $10^{-3}$ & $1.4\times10^{-3}$ & 0.4949 \\
 & Hybrid approach (PCE) & $5\times10^{-4}$ & $9.6\times10^{-4}$ & 0.5316\Bstrut\\\hline\Tstrut
Robust & -- & -- & $1.51\times10^{-2}$ & 0.3380 \\ [1ex] 
\hline 
\end{tabular}
\label{table:ex3_Pf} 
\end{table}

\subsubsection{Case (b): Three-dimensional Beam}

\begin{figure}[!htb]
    \centering
\includegraphics[scale=1]{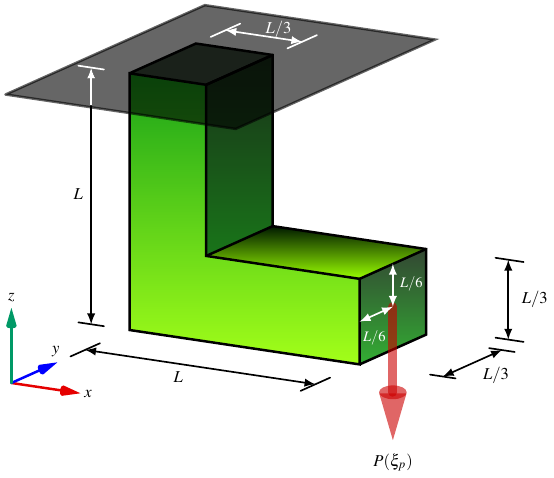}
    \caption{Three dimensional L-shaped beam with uncertain material property and subjected to an uncertain load $P({\color{black}\xi_p})$ at the middle of the right face used in Example III (b). }
    \label{fig:ex3_schem3D}
\end{figure}

\begin{figure}[!htb]
    \centering
    \begin{subfigure}[t]{\textwidth}
    \centering
    \includegraphics[scale=1]{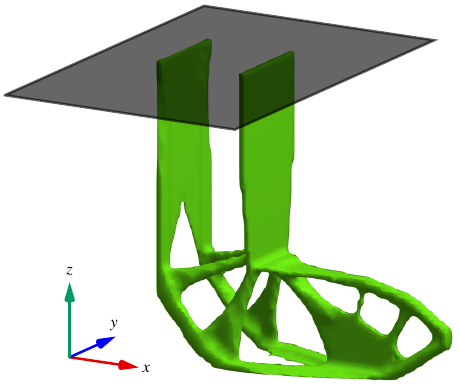}
    
    \caption{RBTO design obtained using the proposed approach} \label{fig:ex3_rbdo_3D} 
    \end{subfigure} \\~\\ 
    \begin{subfigure}[t]{\textwidth}
    \centering
    \includegraphics[scale=1]{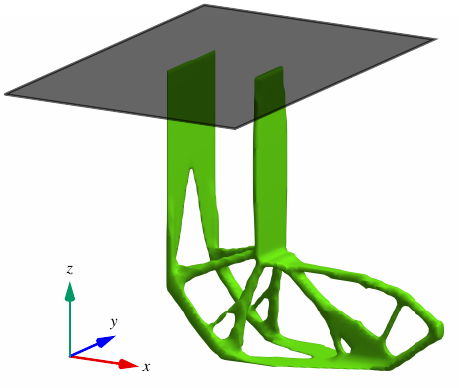}
    
    \caption{Robust design that does not consider the reliability constraint}
    \label{fig:ex3_rto_3D}
    \end{subfigure}
    \caption{Comparison of the RBTO design {\color{black}with a density threshold of 0.7} using a hybrid approach for failure probability estimation for an allowable failure probability $p_a=10^{-2}$ and the robust design {\color{black}with a density threshold of 0.7} that does not include the reliability constraint. } 
    \label{fig:ex3_designs_3D}
\end{figure}

\begin{figure}[!htb]
    \centering
    \begin{subfigure}[t]{\textwidth}
    \centering
    {\includegraphics[scale=0.3]{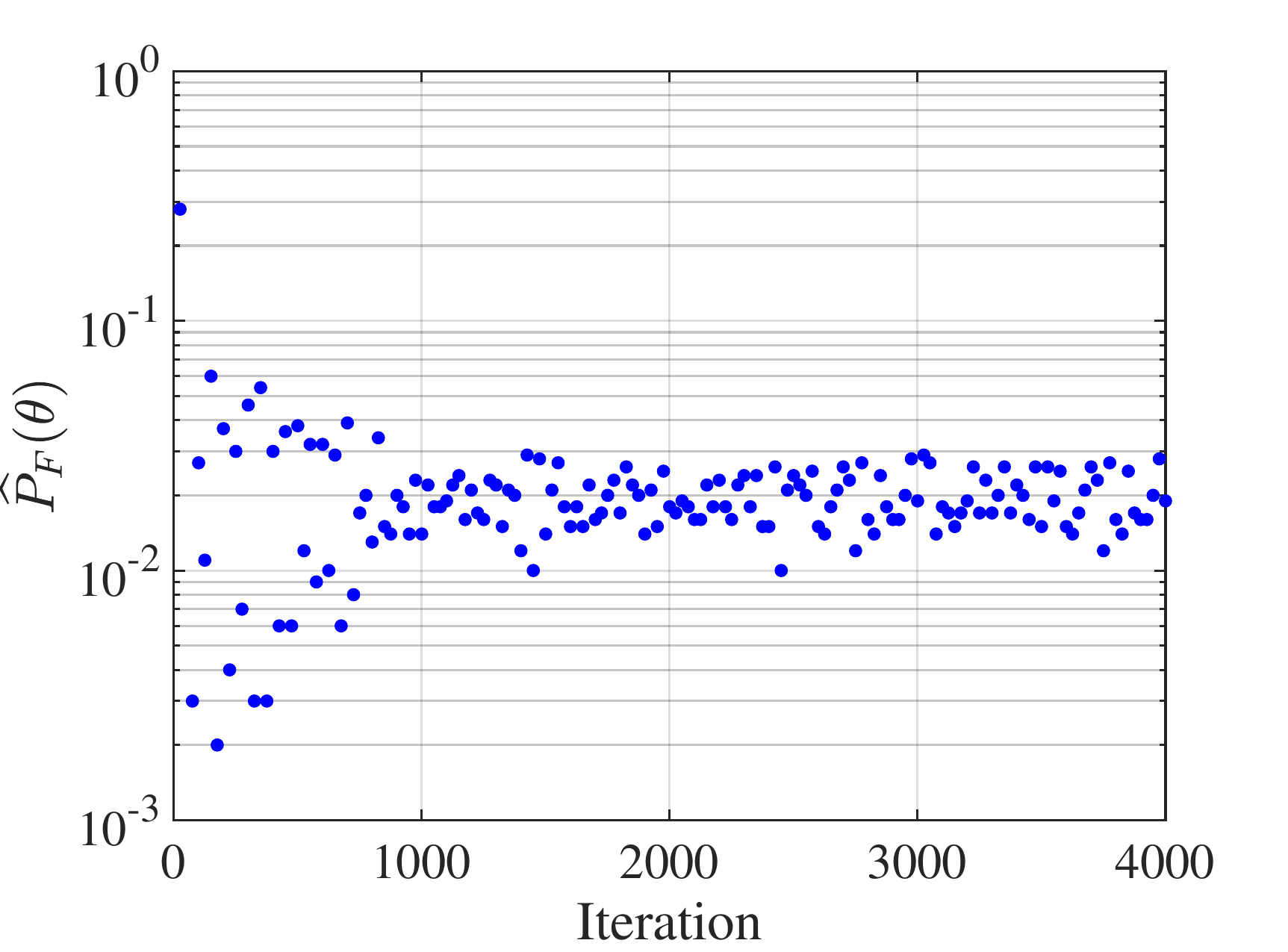}}
    \caption{ The failure probabilities of the design $\widehat{P}_F(\thetaa)$ as estimated during the optimization process} \label{fig:ex3_PfHist_3D} 
    \end{subfigure} \\
    \begin{subfigure}[t]{\textwidth}
    \centering
    {\includegraphics[scale=0.3]{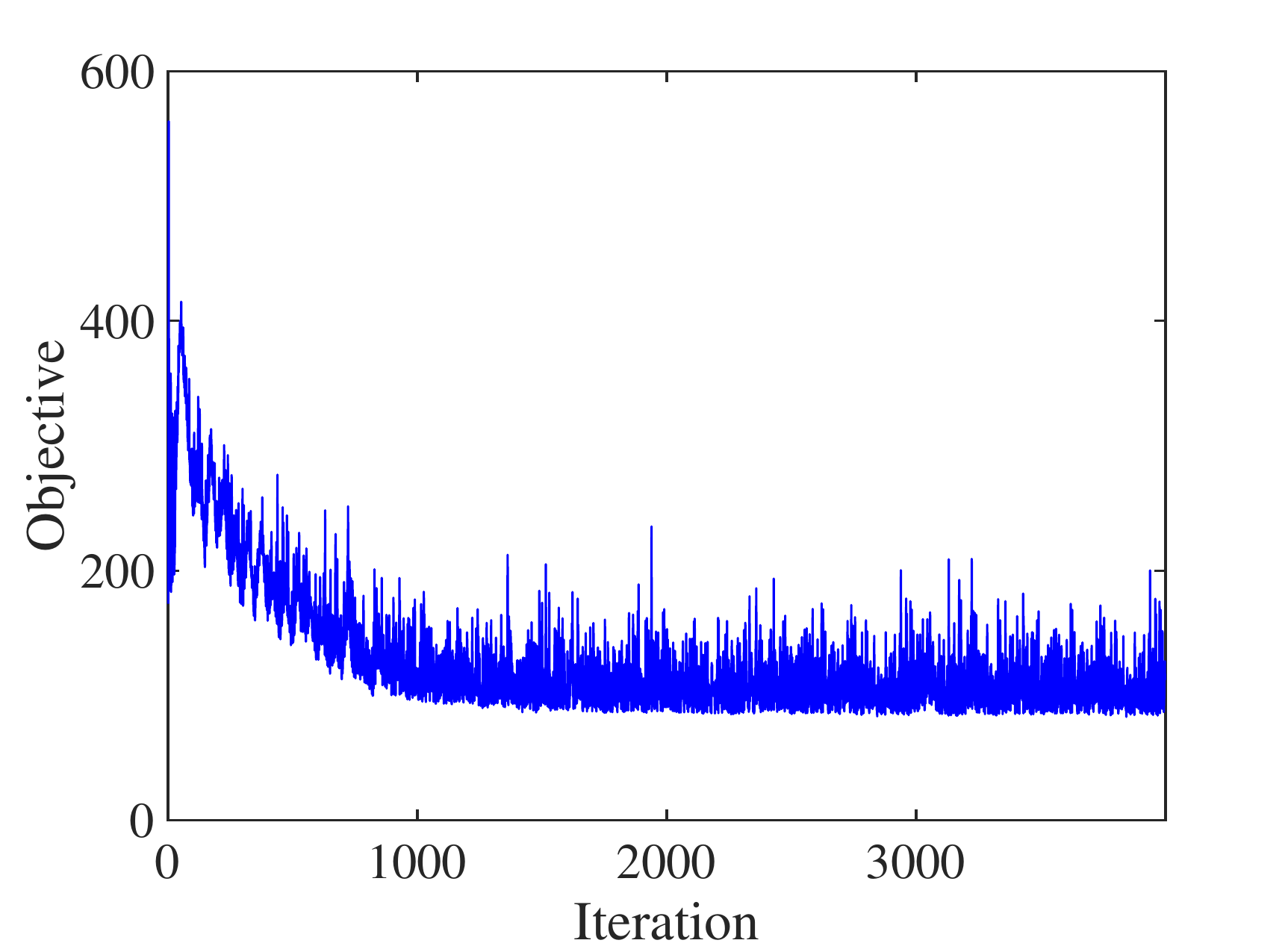}}
    \caption{Objective values for $p_a=10^{-2}$ in the reliability constraint} \label{fig:ex3_ObjHist_3D} 
    \end{subfigure} 
    \caption{Failure probabilities and objective values estimated during the optimization process using a hybrid approach with PCE as surrogate model in Example III (b). } 
    \label{fig:ex3_results_3D}
\end{figure}

Finally, we apply our proposed method for RBTO to the design of the three-dimensional analog of case (a) (see Figure \ref{fig:ex3_schem3D}). Here, we select $\tau=0.01$. The failure event is the compliance value exceeding $C_{\max}=75$, and the maximum allowable probability of failure is $p_a=10^{-2}$. The design domain is discretized into a total of $\sim 70,000$ trilinear elements. 
The design optimization is performed with a step size $\eta = 0.5$ to update the parameters $\thetaa$ for 4000 iterations. We use the same initialization for the design parameters $\thetaa$ as used in Example II. Only one random sample of the uncertain parameters is used in each iteration of the stochastic gradient descent algorithm, i.e., $n=1$. 
We initialize the parameters $\alpha$ and $\betaa$ of $p(\thetaa|F)$ to $10^{-7}$ each,  and set the step size $\eta_F=10^{-7}$. Note that the update in \eqref{eq:ab} uses $w_0=1$ and $w_l=0,l=1,\dots,n_{\thetaa}$ for the first 3000 iterations, and only for the final 1000 iterations, we use $w_0=w_1=\dots=w_{n_{\thetaa}}$ with $\sum_{i=0}^{n_{\thetaa}}w_l=1$. 
The failure probabilities are estimated at every $m=25$ iterations using the hybrid approach with tolerance limit $\gamma=2.5$ and $N=10^3$ in Section \ref{sec:hybrid} (see Algorithm \ref{alg:Hybrid}). The surrogate model is a third order PCE expansion. The coefficients of PCE are estimated using least squares with 16 evaluations of the exact limit state function. In this example, we only use the hybrid approach as the previous examples showed it to be the most economic option compared to the other two approaches. The unconstrained formulation of the optimization problem uses a penalty parameter $\kappa_F=10^5$. Figure \ref{fig:ex3_rbdo_3D} shows the RBTO design obtained from the proposed approach {\color{black}with a density threshold of 0.7}. In contrast to this design, Figure \ref{fig:ex3_rto_3D} shows a robust design {\color{black}with a density threshold of 0.7}, which is obtained by solving the same optimization problem as \eqref{eq:top_obj} but without the reliability constraint. In comparing these two designs, we observe that the RBTO design does not have some of the thin members of the robust design and chooses to use webs instead of circular members in some places. Figure \ref{fig:ex3_results_3D} depicts the estimates of the failure probability and objective during the optimization process, which shows that the results converge after 1500 iterations. {\color{black}Note that the robust design has a failure probability of 0.1048, significantly higher than the allowable limit, when evaluated using the hybrid approach with tolerance limit $\gamma = 2.5$, $N=10^4$, and the PCE surrogate model as before. }
As this case is just a three dimensional extension of case (a), one may expect the RBTO designs obtained here will be very similar to the two-dimensional ones. However, the three dimensional RBTO design has fewer bars than two dimensional RBTO designs in Figure \ref{fig:ex3_results}. Instead, thicker bars are preferred by the optimizer in the three dimensional case. This confirms that a simple extension from two dimension to three dimension is not possible for RBTO designs in this example. However, in three dimensional problems, the evaluation of the limit state function, hence the failure probability, becomes more computationally expensive. The proposed stochastic gradient-based method provides a considerably less expensive option for such cases. 

\section{Conclusions} 

A computational bottleneck in RBTO of large scale structures is the calculation of the statistics, namely the mean of the objective function and constraints, as well as the failure probability. To reduce the computational cost, approximate reliability analysis has been used in the past. In this paper, we propose an alternative approach which is based on a stochastic gradient descent method that updates the design parameters using only $\mathcal{O}(1)$ randomly generated realizations of the objective and constraints at every iteration. This is in contrast to other techniques requiring accurate estimation of the statistics involved using large ensembles of such realizations. Further, we apply Bayes' theorem to a local approximation of the failure probability for a given instance of the design parameters. We assume a parametric exponential form for the probability density function of the design variables within the failure region, {\color{black}and extend it to RBTO problems, where the number of design parameters is very large.} The parameters of this function are updated using stochastic gradient descent as well using $\mathcal{O}(1)$ samples of the design variables from the failure region.

We illustrate the proposed approach using three numerical examples. The first example uses a benchmark problem for a two-bar truss from \cite{rozvany2011analytical}. This example shows that the proposed approach can achieve a reliable design with a fraction of the computational cost of the standard Monte Carlo method. In our second example, we optimize the topology of a beam in the presence of uncertainty in the load and material property. This example shows that a robust design without considering the reliability constraint may have a large probability of failure. The RBTO design adds more features to the design to reduce the chances of failure. The third example uses L-shaped beams in two and three dimensions subjected to uncertain load and material property. Here the reliability-based approach produces a design with thicker members but again with smaller failure probability compared to a robust design. 
{\color{black}In future studies, the proposed method will be compared to other RBTO approaches, e.g., based on FORM/SORM.} 
Further, reliability-based design of coupled multi-physics systems, structural examples with time-variant reliability, {\color{black}i.e., the limit-state function either depends on time or random processes,} and many uncertain parameters will be considered. 

\section*{Acknowledgment} 
The authors acknowledge the support of Defense Advanced Research Projects Agency (DARPA) under the TRADES
program (agreement HR0011-17-2-0022). The opinions and conclusions presented in this paper are those of
the authors and do not necessarily reflect the views of DARPA. 

\section*{Conflicts of Interest}
On behalf of all authors, the corresponding author states that there is no conflict of interest. 

\section*{Replication of Results}
The proposed Algorithm \ref{alg:proposed} in Section \ref{sec:method} has been implemented in MATLAB and will be uploaded to the GitHub page \href{https://github.com/CU-UQ/TOuU}{https://github.com/CU-UQ/TOuU} once the paper is published. 

\bibliographystyle{apalike} 
\bibliography{references.bib}

\begin{thebibliography}{}

\bibitem[Acar and Haftka, 2007]{acar2007reliability}
Acar, E. and Haftka, R.~T. (2007).
\newblock Reliability-based aircraft structural design pays, even with limited
  statistical data.
\newblock {\em Journal of Aircraft}, 44(3):812--823.

\bibitem[Agarwal et~al., 2007]{agarwal2007inverse}
Agarwal, H., Mozumder, C.~K., Renaud, J.~E., and Watson, L.~T. (2007).
\newblock An inverse-measure-based unilevel architecture for reliability-based
  design optimization.
\newblock {\em Structural and Multidisciplinary Optimization}, 33(3):217--227.

\bibitem[Agarwal and Renaud, 2004]{agarwal2004reliability}
Agarwal, H. and Renaud, J. (2004).
\newblock Reliability based design optimization using response surfaces in
  application to multidisciplinary systems.
\newblock {\em Engineering Optimization}, 36(3):291--311.

\bibitem[Alvarez and Carrasco, 2005]{alvarez2005minimization}
Alvarez, F. and Carrasco, M. (2005).
\newblock Minimization of the expected compliance as an alternative approach to
  multiload truss optimization.
\newblock {\em Structural and Multidisciplinary Optimization}, 29(6):470--476.

\bibitem[Aoues and Chateauneuf, 2010]{aoues2010benchmark}
Aoues, Y. and Chateauneuf, A. (2010).
\newblock Benchmark study of numerical methods for reliability-based design
  optimization.
\newblock {\em Structural and Multidisciplinary Optimization}, 41(2):277--294.

\bibitem[Asadpoure et~al., 2011]{asadpoure2011robust}
Asadpoure, A., Tootkaboni, M., and Guest, J.~K. (2011).
\newblock Robust topology optimization of structures with uncertainties in
  stiffness -- {A}pplication to truss structures.
\newblock {\em Computers \& Structures}, 89(11-12):1131--1141.

\bibitem[Au, 2005]{au2005reliability}
Au, S. (2005).
\newblock Reliability-based design sensitivity by efficient simulation.
\newblock {\em Computers \& Structures}, 83(14):1048--1061.

\bibitem[Au and Beck, 2001]{au2001estimation}
Au, S.-K. and Beck, J.~L. (2001).
\newblock Estimation of small failure probabilities in high dimensions by
  subset simulation.
\newblock {\em Probabilistic Engineering Mechanics}, 16(4):263--277.

\bibitem[Bae and Wang, 2002]{bae2002reliability}
Bae, K.-r. and Wang, S. (2002).
\newblock Reliability-based topology optimization.
\newblock In {\em 9th AIAA/ISSMO symposium on multidisciplinary analysis and
  optimization}, page 5542.

\bibitem[Barron and Sheu, 1991]{barron1991approximation}
Barron, A.~R. and Sheu, C.-H. (1991).
\newblock Approximation of density functions by sequences of exponential
  families.
\newblock {\em The Annals of Statistics}, pages 1347--1369.

\bibitem[Basudhar and Missoum, 2008]{basudhar2008adaptive}
Basudhar, A. and Missoum, S. (2008).
\newblock Adaptive explicit decision functions for probabilistic design and
  optimization using support vector machines.
\newblock {\em Computers \& Structures}, 86(19-20):1904--1917.

\bibitem[Beck and de~Santana~Gomes, 2012]{beck2012comparison}
Beck, A.~T. and de~Santana~Gomes, W.~J. (2012).
\newblock A comparison of deterministic, reliability-based and risk-based
  structural optimization under uncertainty.
\newblock {\em Probabilistic Engineering Mechanics}, 28:18--29.

\bibitem[Beck et~al., 2015]{beck2015comparison}
Beck, A.~T., Gomes, W.~J., Lopez, R.~H., and Miguel, L.~F. (2015).
\newblock A comparison between robust and risk-based optimization under
  uncertainty.
\newblock {\em Structural and Multidisciplinary Optimization}, 52(3):479--492.

\bibitem[Beck and Zuev, 2017]{beck2015rare}
Beck, J.~L. and Zuev, K.~M. (2017).
\newblock {\em Stochastic Collocation Methods: A Survey}, pages 1075--1100.
\newblock Springer International Publishing, Cham, Switzerland.

\bibitem[Bect et~al., 2017]{bect2017bayesian}
Bect, J., Li, L., and Vazquez, E. (2017).
\newblock Bayesian subset simulation.
\newblock {\em SIAM/ASA Journal on Uncertainty Quantification}, 5(1):762--786.

\bibitem[Bends{\o}e, 1989]{bendsoe1989optimal}
Bends{\o}e, M.~P. (1989).
\newblock Optimal shape design as a material distribution problem.
\newblock {\em Structural Optimization}, 1(4):193--202.

\bibitem[Beyer and Sendhoff, 2007]{beyer2007robust}
Beyer, H.-G. and Sendhoff, B. (2007).
\newblock Robust optimization--a comprehensive survey.
\newblock {\em Computer Methods in Applied Mechanics and Engineering},
  196(33-34):3190--3218.

\bibitem[Bichon et~al., 2008]{bichon2008efficient}
Bichon, B.~J., Eldred, M.~S., Swiler, L.~P., Mahadevan, S., and McFarland,
  J.~M. (2008).
\newblock Efficient global reliability analysis for nonlinear implicit
  performance functions.
\newblock {\em AIAA Journal}, 46(10):2459--2468.

\bibitem[Bottou, 1999]{bottou_1999}
Bottou, L. (1999).
\newblock {\em On-line Learning and Stochastic Approximations}, page 9–42.
\newblock Publications of the Newton Institute. Cambridge University Press.

\bibitem[Chaudhuri et~al., 2019]{chaudhuri2019reusing}
Chaudhuri, A., Marques, A.~N., Lam, R., and Willcox, K.~E. (2019).
\newblock Reusing information for multifidelity active learning in
  reliability-based design optimization.
\newblock In {\em AIAA Scitech 2019 Forum}, page 1222.

\bibitem[Chen and Chen, 2011]{chen2011new}
Chen, S. and Chen, W. (2011).
\newblock A new level-set based approach to shape and topology optimization
  under geometric uncertainty.
\newblock {\em Structural and Multidisciplinary Optimization}, 44(1):1--18.

\bibitem[Chen et~al., 2010]{chen2010rtso}
Chen, S., Chen, W., and Lee, S. (2010).
\newblock Level set based robust shape and topology optimization under random
  field uncertainties.
\newblock {\em Structural and Multidisciplinary Optimization}, 41(4):507--524.

\bibitem[Cheng et~al., 2006]{cheng2006sequential}
Cheng, G., Xu, L., and Jiang, L. (2006).
\newblock A sequential approximate programming strategy for reliability-based
  structural optimization.
\newblock {\em Computers \& Structures}, 84(21):1353--1367.

\bibitem[Ching and Hsieh, 2007a]{ching2007approximate}
Ching, J. and Hsieh, Y.-H. (2007a).
\newblock Approximate reliability-based optimization using a three-step
  approach based on subset simulation.
\newblock {\em Journal of Engineering Mechanics}, 133(4):481--493.

\bibitem[Ching and Hsieh, 2007b]{ching2007local}
Ching, J. and Hsieh, Y.-H. (2007b).
\newblock Local estimation of failure probability function and its confidence
  interval with maximum entropy principle.
\newblock {\em Probabilistic Engineering Mechanics}, 22(1):39--49.

\bibitem[da~Silva et~al., 2020]{da2020comparison}
da~Silva, G.~A., Cardoso, E.~L., and Beck, A.~T. (2020).
\newblock Comparison of robust, reliability-based and non-probabilistic
  topology optimization under uncertain loads and stress constraints.
\newblock {\em Probabilistic Engineering Mechanics}, 59:103039.

\bibitem[De et~al., 2020a]{de2020topology}
De, S., Hampton, J., Maute, K., and Doostan, A. (2020a).
\newblock Topology optimization under uncertainty using a stochastic
  gradient-based approach.
\newblock {\em Structural and Multidisciplinary Optimization},
  62(5):2255--2278.

\bibitem[De et~al., 2020b]{de2020bi}
De, S., Maute, K., and Doostan, A. (2020b).
\newblock Bi-fidelity stochastic gradient descent for structural optimization
  under uncertainty.
\newblock {\em Computational Mechanics}, 66(4):745--771.

\bibitem[De et~al., 2017]{de2017efficient}
De, S., Wojtkiewicz, S.~F., and Johnson, E.~A. (2017).
\newblock Efficient optimal design and design-under-uncertainty of passive
  control devices with application to a cable-stayed bridge.
\newblock {\em Structural Control and Health Monitoring}, 24(2):e1846.

\bibitem[Dematteis et~al., 2019]{dematteis2019extreme}
Dematteis, G., Grafke, T., and Vanden-Eijnden, E. (2019).
\newblock Extreme event quantification in dynamical systems with random
  components.
\newblock {\em SIAM/ASA Journal on Uncertainty Quantification},
  7(3):1029--1059.

\bibitem[Diwekar, 2020]{Diwekar2020}
Diwekar, U.~M. (2020).
\newblock {\em Optimization Under Uncertainty}, pages 151--215.
\newblock Springer International Publishing, Cham, Switzerland.

\bibitem[Doostan and Owhadi, 2011]{doostan2011non}
Doostan, A. and Owhadi, H. (2011).
\newblock A non-adapted sparse approximation of pdes with stochastic inputs.
\newblock {\em Journal of Computational Physics}, 230(8):3015--3034.

\bibitem[dos Santos et~al., 2018]{dos2018reliability}
dos Santos, R.~B., Torii, A.~J., and Novotny, A.~A. (2018).
\newblock Reliability-based topology optimization of structures under stress
  constraints.
\newblock {\em International Journal for Numerical Methods in Engineering},
  114(6):660--674.

\bibitem[Du and Chen, 2004]{du2004sequential}
Du, X. and Chen, W. (2004).
\newblock Sequential optimization and reliability assessment method for
  efficient probabilistic design.
\newblock {\em Journal of Mechanical Design}, 126(2):225--233.

\bibitem[Dunning and Kim, 2013]{dunning2013robust}
Dunning, P.~D. and Kim, H.~A. (2013).
\newblock Robust topology optimization: minimization of expected and variance
  of compliance.
\newblock {\em AIAA Journal}, 51(11):2656--2664.

\bibitem[Enevoldsen and S{\o}rensen, 1994]{enevoldsen1994reliability}
Enevoldsen, I. and S{\o}rensen, J.~D. (1994).
\newblock Reliability-based optimization in structural engineering.
\newblock {\em Structural Safety}, 15(3):169--196.

\bibitem[Foschi et~al., 2002]{foschi2002reliability}
Foschi, R., Li, H., and Zhang, J. (2002).
\newblock Reliability and performance-based design: a computational approach
  and applications.
\newblock {\em Structural Safety}, 24(2-4):205--218.

\bibitem[Frangopol and Maute, 2003]{frangopol2003life}
Frangopol, D.~M. and Maute, K. (2003).
\newblock Life-cycle reliability-based optimization of civil and aerospace
  structures.
\newblock {\em Computers \& Structures}, 81(7):397--410.

\bibitem[Gano et~al., 2006]{gano2006reliability}
Gano, S.~E., Renaud, J.~E., Agarwal, H., and Tovar, A. (2006).
\newblock Reliability-based design using variable-fidelity optimization.
\newblock {\em Structures and Infrastructure Engineering}, 2(3-4):247--260.

\bibitem[Gasser and Schu{\"e}ller, 1997]{gasser1997reliability}
Gasser, M. and Schu{\"e}ller, G.~I. (1997).
\newblock Reliability-based optimization of structural systems.
\newblock {\em Mathematical Methods of Operations Research}, 46(3):287--307.

\bibitem[Ghanem and Spanos, 2003]{ghanem2003stochastic}
Ghanem, R.~G. and Spanos, P.~D. (2003).
\newblock {\em Stochastic finite elements: a spectral approach}.
\newblock Courier Corporation.

\bibitem[Guest and Igusa, 2008]{guest2008structural}
Guest, J.~K. and Igusa, T. (2008).
\newblock Structural optimization under uncertain loads and nodal locations.
\newblock {\em Computer Methods in Applied Mechanics and Engineering},
  198(1):116--124.

\bibitem[Hadigol and Doostan, 2018]{hadigol2018least}
Hadigol, M. and Doostan, A. (2018).
\newblock Least squares polynomial chaos expansion: A review of sampling
  strategies.
\newblock {\em Computer Methods in Applied Mechanics and Engineering},
  332:382--407.

\bibitem[Haldar and Mahadevan, 2000]{haldar2000reliability}
Haldar, A. and Mahadevan, S. (2000).
\newblock {\em Reliability assessment using stochastic finite element
  analysis}.
\newblock John Wiley \& Sons.

\bibitem[Hasofer and Lind, 1974]{hasofer1974exact}
Hasofer, A.~M. and Lind, N.~C. (1974).
\newblock Exact and invariant second-moment code format.
\newblock {\em Journal of the Engineering Mechanics Division}, 100(1):111--121.

\bibitem[Jalalpour and Tootkaboni, 2016]{jalalpour2016efficient}
Jalalpour, M. and Tootkaboni, M. (2016).
\newblock An efficient approach to reliability-based topology optimization for
  continua under material uncertainty.
\newblock {\em Structural and Multidisciplinary Optimization}, 53(4):759--772.

\bibitem[Jensen, 2005]{jensen2005design}
Jensen, H.~A. (2005).
\newblock Design and sensitivity analysis of dynamical systems subjected to
  stochastic loading.
\newblock {\em Computers \& Structures}, 83(14):1062--1075.

\bibitem[Jensen and Catalan, 2007]{jensen2007effects}
Jensen, H.~A. and Catalan, M.~A. (2007).
\newblock On the effects of non-linear elements in the reliability-based
  optimal design of stochastic dynamical systems.
\newblock {\em International Journal of Non-Linear Mechanics}, 42(5):802--816.

\bibitem[Jung and Cho, 2004]{jung2004reliability}
Jung, H.-S. and Cho, S. (2004).
\newblock Reliability-based topology optimization of geometrically nonlinear
  structures with loading and material uncertainties.
\newblock {\em Finite Elements in Analysis and Design}, 41(3):311--331.

\bibitem[Kale and Haftka, 2008]{kale2008tradeoff}
Kale, A.~A. and Haftka, R.~T. (2008).
\newblock Tradeoff of weight and inspection cost in reliability-based
  structural optimization.
\newblock {\em Journal of Aircraft}, 45(1):77--85.

\bibitem[Kang and Liu, 2018]{kang2018reliability}
Kang, Z. and Liu, P. (2018).
\newblock Reliability-based topology optimization against geometric
  imperfections with random threshold model.
\newblock {\em International Journal for Numerical Methods in Engineering},
  115(1):99--116.

\bibitem[Keshavarzzadeh et~al., 2017]{keshavarzzadeh2017topology}
Keshavarzzadeh, V., Fernandez, F., and Tortorelli, D.~A. (2017).
\newblock Topology optimization under uncertainty via non-intrusive polynomial
  chaos expansion.
\newblock {\em Computer Methods in Applied Mechanics and Engineering},
  318:120--147.

\bibitem[Kharmanda et~al., 2002]{kharmanda2002efficient}
Kharmanda, G., Mohamed, A., and Lemaire, M. (2002).
\newblock Efficient reliability-based design optimization using a hybrid space
  with application to finite element analysis.
\newblock {\em Structural and Multidisciplinary Optimization}, 24(3):233--245.

\bibitem[Kharmanda et~al., 2004]{kharmanda2004reliability}
Kharmanda, G., Olhoff, N., Mohamed, A., and Lemaire, M. (2004).
\newblock Reliability-based topology optimization.
\newblock {\em Structural and Multidisciplinary Optimization}, 26(5):295--307.

\bibitem[Kim et~al., 2006]{kim2006reliability}
Kim, C., Wang, S., Rae, K.-r., Moon, H., and Choi, K.~K. (2006).
\newblock Reliability-based topology optimization with uncertainties.
\newblock {\em Journal of Mechanical Science and Technology}, 20(4):494.

\bibitem[Kouri, 2014]{kouri2014multilevel}
Kouri, D.~P. (2014).
\newblock A multilevel stochastic collocation algorithm for optimization of
  {PDE}s with uncertain coefficients.
\newblock {\em SIAM/ASA Journal on Uncertainty Quantification}, 2(1):55--81.

\bibitem[Kouri et~al., 2013]{kouri2013trust}
Kouri, D.~P., Heinkenschloss, M., Ridzal, D., and van Bloemen~Waanders, B.~G.
  (2013).
\newblock A trust-region algorithm with adaptive stochastic collocation for
  {PDE} optimization under uncertainty.
\newblock {\em SIAM Journal on Scientific Computing}, 35(4):A1847--A1879.

\bibitem[Kuschel and Rackwitz, 1997]{kuschel1997two}
Kuschel, N. and Rackwitz, R. (1997).
\newblock Two basic problems in reliability-based structural optimization.
\newblock {\em Mathematical Methods of Operations Research}, 46(3):309--333.

\bibitem[Li and Xiu, 2010]{li2010evaluation}
Li, J. and Xiu, D. (2010).
\newblock Evaluation of failure probability via surrogate models.
\newblock {\em Journal of Computational Physics}, 229(23):8966--8980.

\bibitem[Li and Zhang, 2020]{li2020momentum}
Li, W. and Zhang, X.~S. (2020).
\newblock Momentum-based accelerated mirror descent stochastic approximation
  for robust topology optimization under stochastic loads.
\newblock {\em arXiv preprint arXiv:2008.13284}.

\bibitem[Lopez and Beck, 2012]{lopez2012reliability}
Lopez, R.~H. and Beck, A.~T. (2012).
\newblock Reliability-based design optimization strategies based on form: a
  review.
\newblock {\em Journal of the Brazilian Society of Mechanical Sciences and
  Engineering}, 34(4):506--514.

\bibitem[Luenberger and Ye, 1984]{luenberger1984linear}
Luenberger, D.~G. and Ye, Y. (1984).
\newblock {\em Linear and nonlinear programming}, volume~2.
\newblock Springer.

\bibitem[Luo et~al., 2014]{luo2014reliability}
Luo, Y., Zhou, M., Wang, M.~Y., and Deng, Z. (2014).
\newblock Reliability based topology optimization for continuum structures with
  local failure constraints.
\newblock {\em Computers \& Structures}, 143:73--84.

\bibitem[Madsen et~al., 2006]{madsen2006methods}
Madsen, H.~O., Krenk, S., and Lind, N.~C. (2006).
\newblock {\em Methods of structural safety}.
\newblock Courier Corporation.

\bibitem[Maute, 2014]{maute2014touu}
Maute, K. (2014).
\newblock Topology optimization under uncertainty.
\newblock In {\em Topology Optimization in Structural and Continuum Mechanics},
  pages 457--471. Springer.

\bibitem[Maute and Frangopol, 2003]{maute2003reliability}
Maute, K. and Frangopol, D.~M. (2003).
\newblock Reliability-based design of {MEMS} mechanisms by topology
  optimization.
\newblock {\em Computers \& Structures}, 81(8-11):813--824.

\bibitem[Melchers and Beck, 2018]{melchers2018structural}
Melchers, R.~E. and Beck, A.~T. (2018).
\newblock {\em Structural reliability analysis and prediction}.
\newblock John Wiley \& Sons.

\bibitem[Meng and Keshtegar, 2019]{meng2019adaptive}
Meng, Z. and Keshtegar, B. (2019).
\newblock Adaptive conjugate single-loop method for efficient reliability-based
  design and topology optimization.
\newblock {\em Computer Methods in Applied Mechanics and Engineering},
  344:95--119.

\bibitem[Missoum et~al., 2007]{missoum2007convex}
Missoum, S., Ramu, P., and Haftka, R.~T. (2007).
\newblock A convex hull approach for the reliability-based design optimization
  of nonlinear transient dynamic problems.
\newblock {\em Computer Methods in Applied Mechanics and Engineering},
  196(29-30):2895--2906.

\bibitem[Moustapha and Sudret, 2019]{moustapha2019surrogate}
Moustapha, M. and Sudret, B. (2019).
\newblock Surrogate-assisted reliability-based design optimization: a survey
  and a unified modular framework.
\newblock {\em Structural and Multidisciplinary Optimization}, pages 1--20.

\bibitem[Nguyen et~al., 2011]{nguyen2011single}
Nguyen, T.~H., Song, J., and Paulino, G.~H. (2011).
\newblock Single-loop system reliability-based topology optimization
  considering statistical dependence between limit-states.
\newblock {\em Structural and Multidisciplinary Optimization}, 44(5):593--611.

\bibitem[Nikolaidis et~al., 2004]{nikolaidis2004comparison}
Nikolaidis, E., Chen, S., Cudney, H., Haftka, R.~T., and Rosca, R. (2004).
\newblock Comparison of probability and possibility for design against
  catastrophic failure under uncertainty.
\newblock {\em Journal of Mechanical Design}, 126(3):386--394.

\bibitem[Ramu et~al., 2006]{ramu2006inverse}
Ramu, P., Qu, X., Youn, B.~D., Haftka, R.~T., and Choi, K.~K. (2006).
\newblock Inverse reliability measures and reliability-based design
  optimisation.
\newblock {\em International Journal of Reliability and Safety},
  1(1-2):187--205.

\bibitem[Rozvany and Maute, 2011]{rozvany2011analytical}
Rozvany, G.~I. and Maute, K. (2011).
\newblock Analytical and numerical solutions for a reliability-based benchmark
  example.
\newblock {\em Structural and Multidisciplinary Optimization}, 43(6):745--753.

\bibitem[Sigmund, 2001]{sigmund200199}
Sigmund, O. (2001).
\newblock A 99 line topology optimization code written in {Matlab}.
\newblock {\em Structural and Multidisciplinary Optimization}, 21(2):120--127.

\bibitem[Sigmund and Maute, 2013]{sigmund2013topology}
Sigmund, O. and Maute, K. (2013).
\newblock Topology optimization approaches.
\newblock {\em Structural and Multidisciplinary Optimization},
  48(6):1031--1055.

\bibitem[Silva et~al., 2010]{silva2010component}
Silva, M., Tortorelli, D.~A., Norato, J.~A., Ha, C., and Bae, H.-R. (2010).
\newblock Component and system reliability-based topology optimization using a
  single-loop method.
\newblock {\em Structural and Multidisciplinary Optimization}, 41(1):87--106.

\bibitem[Suryawanshi and Ghosh, 2016]{suryawanshi2016reliability}
Suryawanshi, A. and Ghosh, D. (2016).
\newblock Reliability based optimization in aeroelastic stability problems
  using polynomial chaos based metamodels.
\newblock {\em Structural and Multidisciplinary Optimization},
  53(5):1069--1080.

\bibitem[Taflanidis and Beck, 2008a]{taflanidis2008stochastic}
Taflanidis, A. and Beck, J. (2008a).
\newblock Stochastic subset optimization for optimal reliability problems.
\newblock {\em Probabilistic Engineering Mechanics}, 23(2-3):324--338.

\bibitem[Taflanidis and Beck, 2008b]{taflanidis2008efficient}
Taflanidis, A.~A. and Beck, J.~L. (2008b).
\newblock An efficient framework for optimal robust stochastic system design
  using stochastic simulation.
\newblock {\em Computer Methods in Applied Mechanics and Engineering},
  198(1):88--101.

\bibitem[Tong et~al., 2020a]{tong2020optimization}
Tong, S., Subramanyam, A., and Rao, V. (2020a).
\newblock Optimization under rare chance constraints.
\newblock {\em arXiv preprint arXiv:2011.06052}.

\bibitem[Tong et~al., 2020b]{tong2020extreme}
Tong, S., Vanden-Eijnden, E., and Stadler, G. (2020b).
\newblock Extreme event probability estimation using pde-constrained
  optimization and large deviation theory, with application to tsunamis.
\newblock {\em arXiv preprint arXiv:2007.13930}.

\bibitem[Tootkaboni et~al., 2012]{tootkaboni2012topology}
Tootkaboni, M., Asadpoure, A., and Guest, J.~K. (2012).
\newblock Topology optimization of continuum structures under uncertainty -- {A
  Polynomial Chaos} approach.
\newblock {\em Computer Methods in Applied Mechanics and Engineering},
  201:263--275.

\bibitem[Torii et~al., 2016]{torii2016robust}
Torii, A.~J., Novotny, A.~A., and dos Santos, R.~B. (2016).
\newblock Robust compliance topology optimization based on the topological
  derivative concept.
\newblock {\em International Journal for Numerical Methods in Engineering},
  106(11):889--903.

\bibitem[Tu et~al., 1999]{tu1999new}
Tu, J., Choi, K., and Park, Y. (1999).
\newblock A new study on reliability-based design optimization.
\newblock {\em Journal of Mechanical Design}, 121(4):557--564.

\bibitem[Valdebenito and Schu{\"e}ller, 2010]{valdebenito2010survey}
Valdebenito, M.~A. and Schu{\"e}ller, G.~I. (2010).
\newblock A survey on approaches for reliability-based optimization.
\newblock {\em Structural and Multidisciplinary Optimization}, 42(5):645--663.

\bibitem[Xiu and Karniadakis, 2002]{xiu2002wiener}
Xiu, D. and Karniadakis, G.~E. (2002).
\newblock The {W}iener--{A}skey polynomial chaos for stochastic differential
  equations.
\newblock {\em SIAM Journal on Scientific Computing}, 24(2):619--644.

\bibitem[Yang and Gu, 2004]{yang2004experience}
Yang, R. and Gu, L. (2004).
\newblock Experience with approximate reliability-based optimization methods.
\newblock {\em Structural and Multidisciplinary Optimization}, 26(1):152--159.

\bibitem[Zhang and Foschi, 2004]{zhang2004performance}
Zhang, J. and Foschi, R.~O. (2004).
\newblock Performance-based design and seismic reliability analysis using
  designed experiments and neural networks.
\newblock {\em Probabilistic Engineering Mechanics}, 19(3):259--267.

\bibitem[Zuev et~al., 2012]{zuev2012bayesian}
Zuev, K.~M., Beck, J.~L., Au, S.-K., and Katafygiotis, L.~S. (2012).
\newblock Bayesian post-processor and other enhancements of subset simulation
  for estimating failure probabilities in high dimensions.
\newblock {\em Computers \& Structures}, 92:283--296.

\end{thebibliography}

\end{document}